\documentclass[twoside,10pt,draft]{book}
\usepackage{amsthm}
\usepackage{amsmath}
\usepackage{amssymb}
\usepackage{latexsym}
\usepackage{enumerate}

\numberwithin{equation}{chapter}

\newtheorem{theorem}{Theorem}[chapter]

\newtheorem{corollary}[theorem]{Corollary}
\newtheorem{proposition}[theorem]{Proposition}
\newtheorem{lemma}[theorem]{Lemma}
\theoremstyle{definition}
\newtheorem{definition}[theorem]{Definition}
\theoremstyle{remark}
\newtheorem{remark}[theorem]{Remark}

\newtheorem{example}[theorem]{Example}
\newtheorem{exercise}[theorem]{Exercise}

\newcommand\cS{\mathcal{S}}
\newcommand\cC{\mathcal{C}}
\newcommand{\cO}{\mathcal{O}}

\newcommand{\cU}{\mathcal{U}}

\newcommand{\R}{\mathbb{R}}
\newcommand{\C}{\mathbb{C}}
\newcommand{\Z}{\mathbb{Z}}
\newcommand{\Q}{\mathbb{Q}}
\newcommand{\bbT}{\mathbb{T}}
\newcommand\lie[1]{\mathfrak{#1}}

\newcommand{\fh}{\lie{h}}
\newcommand{\fg}{\lie{g}}

\def    \inv    {^{-1}}

\newcommand\im{\mathop{\rm im}\nolimits}

\newcommand{\tPhi} {\tilde{\Phi}}

\def	\tU	{{ \tilde{U}  }}
\def	\tV	{{ \tilde{V}  }}     

\newcommand\tcU {{{\mathcal U}}}
\newcommand\tcV {{{\mathcal V}}}
\newcommand\tf {{\tilde{f}}}
\newcommand\tsigma {{\tilde{\sigma}}}

\newcommand	{\printname}[1]	{}

\newcommand{\labell}[1]	{\label{#1}\printname{#1}}

\newcommand{\clearemptydoublepage}{\newpage{\pagestyle{empty}\cleardoublepage}}

\begin{document}

\title{Geodesic flows and contact 	toric manifolds \\
\vspace*{10ex}}
\author{Eugene Lerman \thanks{
Department of Mathematics, University of Illinois, Urbana, IL 61801, USA}
\thanks{Email: {\tt
lerman@math.uiuc.edu} }\thanks{Partially supported by NSF grant DMS - 980305 and R. Kantorovitz.}}
\date{January 2002}

\maketitle

\pagenumbering{roman}

\thispagestyle{plain}

\section*{Forward}
These notes are based on five 1.5 hour lectures on torus actions on
contact manifolds delivered at the summer school on Symplectic
Geometry of Integrable Hamiltonian Systems at Centre de Recerca
Matem\`atica in Barcelona in July 2001. Naturally the notes contain
more material that could have been delivered in 7.5 hours.  I am
grateful to Carlos Curr\`as-Bosch and Eva Miranda, the organizers of
the summer school, for their kind invitation to teach a course.
Thanks are also due to the staff of the CRM without whom the summer
school would not have been a success.\\

The main theme of these notes is the topological study of contact
toric manifolds, a relatively new class of manifolds that I find very
interesting. A motivation for studying these manifolds comes from
completely integrable systems $\{f_1, \ldots, f_n\}$ on punctured
cotangent bundles where each function $f_i$ is homogeneous of degree 1
(one can think of $f_i$'s as symbols of first order
pseudo-differential operators, but this is not essential).  A
punctured cotangent bundle is a symplectic cone whose base is
naturally a contact manifold (this is explained in detail in
Chapter~\ref{sec2}).  This observation leads to studying completely
integrable systems on contact manifolds, whatever those are.

The simplest (symplectic) completely integrable systems are the ones
with global action-angle coordinates.  The next simplest case is that
of Hamiltonian torus actions.  If the phase space is compact one ends
up with (compact) symplectic toric manifolds.  This is the theme of
Ana Cannas's lectures delivered at the summer school.  The
corresponding case in the contact category is that of compact toric
manifolds.

We will use the excuse of studying completely integrable geodesic
flows with homogeneous integrals to introduce various ideas essential
for the classification of contact and symplectic toric manifolds.
More specifically we will discuss in these notes contact moment maps,
slices for group actions, sheaves and \v Cech cohomology, 
orbifolds and Morse theory on orbifolds.

\clearemptydoublepage

\tableofcontents

\clearemptydoublepage

\pagestyle{headings}

\pagenumbering{arabic}

\setcounter{page}{1}

\chapter{Introduction}

We start with an innocuous sounding problem.

\noindent
{\bf Problem}\quad
Consider the cotangent bundle of the $n$-torus minus the zero section
$T^*\bbT^n \smallsetminus 0$.  That is, consider the manifold 
$$ 
M =\{(q, p) \in \R^n/\Z^n \times \R^n \mid p \not = 0 \} .  
$$

Suppose further that the torus $G=\bbT^n = \R^n/\Z^n$ acts on
$M$ effectively and preserves the standard symplectic form $\omega =
\sum dp_i \wedge dq_i$ (i.e., the action is symplectic).  Suppose further 
that the action of $G$ commutes with dilations, i.e., the action
$\rho$ of $\R$ on $M$ given by 
$$
\rho _\lambda (q, p) = (q, e^\lambda p).
$$
Is the action of $G$ necessarily free?\\

\begin{remark}\labell{remarks}
\begin{enumerate}

\item
Recall that an action of the group $G$ on a manifold $M$ is {\bf effective}
if the only element of $G$ that fixes all the points of $M$ is the
identity.

\item There is an ``obvious'' action of $\bbT^n$ on $M$ which has the above
properties and is free  --- it is the lift of  left multiplication:
$$
a \cdot (q, p) = (a  q, p),
$$
where $a\in\bbT^n$, $(q, p) \in M$.  
The issue is
whether an {\bf arbitrary} action of $\bbT^n$ which is effective,
symplectic and commutes with dilations is necessarily free.

\item 
We will see later (Proposition~\ref{mmap-on-cones}) that an action of
a Lie group $G$ on a punctured cotangent bundle which is symplectic and
commutes with dilations is necessarily {\bf Hamiltonian}. That is,
there is a {\bf moment map} $\Phi :M \to \fg^*$.  Recall the
definition of the moment map for a symplectic action of a Lie group
$G$ on a symplectic manifold $(M, \omega)$: for any vector $X$ in the
Lie algebra $\fg$ of $G$ 
$$
d \langle \Phi, X\rangle = \omega (X_M , \cdot ), 
$$ 
where $\langle \cdot, \cdot \rangle : \fg^* \times \fg \to \R$ is
the canonical pairing and $X_M$ is the vector field induced by the
infinitesimal action of $X$: $X_M (x) = \left. \frac{d}{dt}\right|_{t
= 0} (\exp tX) \cdot x$.

\end{enumerate}
\end{remark}

The problem is due to John Toth and Steve Zelditch \cite{TZ}.  The
context is classical and quantum integrability of geodesic flows.
Recall that for a manifold $Q$ with a Riemannian metric $g$ the
corresponding {\bf geodesic flow} is the flow on the cotangent bundle
$T^*Q$ of the Hamiltonian vector field $X_h$ of the function $h = h_g
\in C^\infty (T^*Q)$ which is the square root of the energy:
$$ 
h_g (q, p) = (g_q^* (p, p))^{1/2} 
$$ 
for all $q\in Q$, $p\in T^*_q Q$.  Here $g^*$ denotes the inner
product on the cotangent bundle $T^*Q \to Q$ dual to the inner product
$g$; $g^*$ is the so called dual metric.  For an analyst the function
$h$ is the principal symbol of the square root of the Laplace operator
$\sqrt{\Delta}$ defined by the metric $g$.  A precise definition of
$\sqrt{\Delta}$ will play no role in these notes.  For a Riemannian
geometer it's important that the integral curves of the vector field
$X_h$ project down to geodesics on the manifold $Q$.  Toth and
Zelditch were interested in the meaning of $L^\infty$ boundedness of
the $L^2$-normalized eigenfunctions of $\sqrt{\Delta}$.  They observed
that the question is easier if $\sqrt{\Delta}$ is quantum completely
integrable.  Again, it will not be important to us as to what that
means precisely.  What will matter is that quantum integrability of
$\sqrt{\Delta}$ implies (classical) homogeneous complete integrability
of the geodesic flow.  Namely, it implies that there exist functions
$f_1 = h, f_2,
\ldots, f_n \in C^\infty (T^*Q \smallsetminus 0)$, $n = \dim Q$,\footnote{
$T^* Q \smallsetminus 0$ denotes the punctured cotangent bundle of
$Q$, that is, $T^*Q$ with the zero section deleted.}
such that
\begin{enumerate}
\item 
the functions $f_1, \ldots, f_n$ are functionally independent on an
open dense set $U\subset T^*Q\smallsetminus 0$, i.e., $df_1 \wedge
\ldots \wedge df_n
\not = 0$ on $U$;

\item 
the functions Poisson commute with each other: $\{f_i, f_j\} = 0$ for
all $1\leq i, j \leq n$;

\item 
the functions $f_i$ are homogeneous of degree 1:
$$
\rho_\lambda ^*f_i =e ^\lambda f_i
$$ 
for all $\lambda \in \R$, where $\rho_\lambda : T^*Q \smallsetminus
0 \to T^*Q \smallsetminus 0$ again denotes the dilation $\rho_\lambda (q, p) =
(q, e^\lambda p)$, $q\in Q$, $p\in T^*_q Q $.
\end{enumerate}

\begin{exercise}
Suppose $f$ is a smooth function on the punctured cotangent bundle
$T^*Q \smallsetminus 0$ of a manifold $Q$ which is homogeneous of
degree 1, i.e., $\rho_\lambda ^*f =e ^\lambda f$ for all dilations
$\rho_\lambda$.  Show that its Hamiltonian vector field $X_f$ (relative 
to the standard symplectic structure on $T^*Q$) satisfies 
$$ 
d\rho_\lambda (X_f) = X_f \circ \rho_\lambda .  
$$ 
Conclude that the flow of $X_h$ commutes with dilations.
\end{exercise}

Given a completely integrable system, we know that locally around any
generic point there exist action-angle variables.  Toth and Zelditch
observed that things are considerably simpler if the action-angle
variables are global.  Then there exists an effective action of a
torus $\bbT^n =\R^n /\Z^n$ on $T^*Q \smallsetminus 0$ preserving the
function $h$ and the symplectic form, and commuting with dilations
$\rho_\lambda$.  Toth and Zelditch proved ({\em op.\, cit.}):

\begin{theorem}
Suppose that the Lie group $G=\bbT^n$ acts effectively
on $M = T^* \bbT^n \smallsetminus 0$, preserving the standard
symplectic form and the function $h_g (q, p) = (g^*_q (p,p))^{1/2}$
for some metric $g$ on $\bbT^n$.  Suppose further the action commutes
with dilations.  If the action of $G$ is free then the metric $g$ is
flat, that is, 
$$ 
g = \sum g_{ij}\, dq_i \otimes dq _j 
$$ 
for some {\em constants }$g_{ij}$ (with $g_{ij} = g_{ji}$).
\end{theorem}
\noindent
The eigenfunctions of a flat metric Laplace operator on a torus are
well understood.

Let us now go back to the problem.  The answer to the question is yes
\cite{LS}:

\begin{theorem}\labell{main}
Suppose the Lie group $G=\bbT^n =
\R^n /\Z^n$ acts effectively on the punctured cotangent bundle $M =
T^* \bbT^n \smallsetminus 0$ preserving the standard symplectic form
and commuting with dilations $\rho_\lambda : M \to M$, $\rho_\lambda
(q, p) = (q, e^\lambda p)$.  Then the action of $G$ is free.
\end{theorem}

The main purpose of these notes is to explain why Theorem~\ref{main}
is true.  I will now try to motivate the proof and to put it in a
broader context.  As was remarked previously
(Remark~\ref{remarks}(3)), the action of the torus $G$ on $M = T^*
\bbT^n \smallsetminus 0$ is Hamiltonian.  By the dimension count the
manifold $M$ together with its natural symplectic structure and the
action of $G$ is a symplectic toric manifold.  Recall the definition.

\begin{definition}
A {\bf symplectic toric manifold} is a triple $(M, \omega, \Phi :M \to
\fg^*)$ where $M$ is a manifold, $\omega$ is a symplectic form on $M$,
and $\Phi $ is a moment map for an effective Hamiltonian action of a
torus $G$ on $(M, \omega)$ satisfying $2 \dim G = \dim M$.
\end{definition}

{\em Compact} symplectic toric manifolds are well understood thanks to
a classification theorem of Delzant \cite{D}, which says that all such
manifolds are classified by the images of the corresponding moment
maps.  Note that by the Atiyah - Guillemin - Sternberg convexity
theorem \cite{A, GS}, the images are convex rational polytopes.
Delzant proved that in the toric case the polytopes are
simple\footnote{ A polytope in an $n$-dimensional real vector space is
{\bf simple} if there are exactly $n$ edges meeting at each vertex.
Equivalently, all the supporting hyperplanes are in general position.
Thus a cube and a tetrahedron are simple and an octahedron is not.  }
and additionally satisfy certain integrality conditions. Finally any
simple polytope satisfying the integrality conditions occurs as an
image of the moment map for a compact symplectic toric manifold.

\begin{exercise}
Show that a Hamiltonian torus action on a {\em compact} symplectic
manifold is never free. Hint: any smooth function on a compact
manifold has a critical point (in fact it has at least two --- a
maximum and a minimum).
\end{exercise}

The manifold $M = T^*\bbT^n \smallsetminus 0$ we are interested in is
not compact.  Worse, we will see that the corresponding moment map
$\Phi : M \to \fg^* $ is homogeneous (Proposition~\ref{mmap-on-cones}):
$$
\Phi (\rho _\lambda (m)) = e^\lambda \Phi (m)
$$ for all $m\in M$, $\lambda \in \R$.  This is a bit of bad news ---
Morse theory is an essential ingredient in the proof of the
Atiyah-Guillemin-Sternberg convexity theorem and hence of Delzant's
classification. Morse functions on a noncompact manifold which are not
bounded either above or below are in practice impossible to work with.

On the other hand, because of homogeneity, the moment map descends to
the quotient of $M$ by the action of $\R$.  The quotient is
diffeomorphic to the co-sphere bundle $S^*\bbT^n : = \{(q, p) \in T^*
\bbT^n \mid g^*_q (p, p) = 1\}$ for some dual metric $g^*$.  Since
$S^* \bbT^n$ is odd-dimensional, it is not a symplectic
manifold.\footnote{The manifold $S^* \bbT^n$ is contact. See next
Chapter and the Appendix.}

It is also easy to see in an example that the map induced on $S^*
\bbT^n$ by a homogeneous moment map on $T^*\bbT^n \smallsetminus 0$
behaves rather strangely if one is used to symplectic moment maps:

\begin{example}
Consider the standard action of $G= \bbT^n$ on the cotangent bundle
$T^* \bbT^n$, the lift of the left multiplication:
$$
a \cdot (q, p ) = (a \cdot q, p),
$$
The corresponding moment map $\Phi : T^* G = G \times \fg^* \to \fg^*$
is given by $\Phi (q, p) = p$.  Fix the standard metric on $G$ and
identify the dual of the Lie algebra $\fg^*$ with $\R^n$.  Then the
co-sphere bundle $S^* G$ is $\bbT^n \times S^{n-1}$.  The map $\Phi' =
\Phi |_{S^*G}: S^* G \to \fg^*$ is also given by
$\Phi' (q, p) = p$.  Note that for any nonzero vector $X\in \fg$ the
function $\langle \Phi', X\rangle$ has exactly two critical manifolds
even though the action of $G$ on $S^*G$ is free!
\end{example}

 Compare this with the symplectic situation where critical points of
 components of moment maps are points with nontrivial isotropy groups.
 What are we dealing with? We are dealing with moment maps for group
 actions on contact manifolds.

We finish the section by sketching a strategy for our proof of
Theorem~\ref{main}.  The effective action of the torus $G$ on $M = T^*
\bbT^n \smallsetminus 0$ descends to an effective action on the
quotient $B =S^* \bbT ^n$ of $M$ by dilations.  The action of $G$ on
$B$ preserves a contact structure $\xi$ (see
Definition~\ref{def-cont-structure} below) making $(B, \xi)$ into a
compact connected contact toric $G$-manifold (c.c.c.t.m., see
Definition~\ref{def-ctm}).  We then study all c.c.c.t.m.'s with a
non-free torus actions  and argue that none of them can
have the homotopy type of $B = \bbT^n
\times S^{n-1}$.

\pagestyle{headings}

\clearemptydoublepage

\chapter{Symplectic cones and  contact manifolds} \label{sec2}

In this section we define symplectic cones, contact forms and contact
structures.  Given a symplectic cone we show how to construct the
corresponding contact manifold, and conversely, given a contact
manifold we construct the corresponding symplectic cone. Thus
symplectic manifolds and contact manifolds are ``the same thing.''
Next we show that a symplectic action of a Lie group on a symplectic
cone induces a contact action on the corresponding contact manifold.
This will give us tools to set up a proof of Theorem~\ref{main} as a
study of contact toric manifolds.  The material in this section is
fairly well known. We now start by defining symplectic cones.

\begin{definition}  
A symplectic manifold $(M, \omega)$ is a {\bf symplectic cone} if
\begin{itemize}
\item the manifold $M$ is a principal $\R$ bundle over some manifold $B$, 
called the {\bf base} of the cone,  and

\item
the action of the real line $\R$ expands the symplectic form
exponentially.  That is, $\rho_\lambda ^*\omega = e^\lambda \omega$,
where $\rho_\lambda$ denotes the diffeomorphism define by $\lambda
\in \R$.
\end{itemize}
\end{definition}
\begin{definition}\labell{def-proper}
Recall that a map $f: X \to Y$ between two topological
spaces is {\bf proper} if the preimage of a compact set under $f$ is
compact.
 An action  of a Lie group $G$ on a manifold $M$ 
is {\bf proper} if  the map
$G\times M \ni (g, m) \mapsto (g\cdot m, m) \in M\times M$ is proper.
\end{definition}

By a theorem of Palais \cite{P} the quotient $M/G$ of a manifold $M$ by a free
proper action of a Lie group $G$ is a manifold and the orbit map $M
\to M/G$ makes $M$ into a principal $G$ bundle (see also 
Remark~\ref{rmrk4.7} below).  It follows that if a 
symplectic manifold $(M, \omega)$ has a complete vector field $X$ with
the following two properties:
\begin{enumerate}
\item
the action of $\R$ induced by the flow of $X$ is proper, and
\item
the Lie derivative of the symplectic form $\omega$ with respect to
the vector field $X$ is again $\omega$: $L_X \omega = \omega$,
\end{enumerate}
then $(M, \omega)$ is a symplectic cone relative to the induced action of $\R$.
This gives us an equivalent definition of a symplectic cone.
\begin{definition}
A {\bf symplectic cone} is a triple $(M, \omega, X)$ where $M$ is a
manifold, $\omega$ is a symplectic form on $M$, $X$ is a vector field
on $M$ generating a proper action of the reals $\R$ such that $L_X
\omega = \omega$.
\end{definition}

\begin{example}
Let $(V, \omega_V)$ be a symplectic vector space.  The manifold $M =
V\smallsetminus \{0\}$ is a symplectic cone with the action of $\R$
given by $\rho_\lambda (v) = e^\lambda v$.  Clearly $\rho_\lambda ^*
\omega_V = e^\lambda \omega_V$. The base is a sphere.
\end{example}

\begin{example}
Let $Q$ be a manifold.  Denote the cotangent bundle of $Q$ with the
zero section deleted by $T^*Q \smallsetminus 0$.  There is a natural
free action of the reals $\R$ on the manifold $M := T^*Q
\smallsetminus 0$ given by dilations 
$\rho_\lambda (q, p) = (q, e^\lambda p)$.  It expands the standard
symplectic form on the cotangent bundle exponentially.  Thus $T^*Q
\smallsetminus 0$ is naturally a symplectic cone.  The base is the 
co-sphere bundle $S^* Q$. 
\end{example}

\begin{proposition}\labell{mmap-on-cones}
Suppose $(M, \omega, X)$ is a symplectic cone and suppose a Lie group
$G$ acts on $M$ preserving the symplectic form $\omega$ and the
expanding vector field $X$.  Then the action of $G$ on the symplectic
manifold $(M, \omega)$ is Hamiltonian.  Moreover we may choose the
moment map $\Phi :M \to \fg^*$ to be homogeneous of degree 1, i.e., $$
\Phi (\rho_\lambda (m)) = e^\lambda \Phi (m)
$$ for all $\lambda \in \R$ and $m\in M$.  Here $\rho_\lambda$ denotes
the action of $\R$ generated by $X$, that is, the time $\lambda$ flow of $X$.
\end{proposition}

\begin{proof}
Note first that since $L_X \omega = \omega$ and since $d\omega =0$, it
follows from Cartan's formula ($L_X \omega = \iota (X) d \omega + d
\iota (X) \omega$) that $d (\iota (X) \omega) = \omega$.  Since the
action of $G$ preserves $X$ and $\omega$, it preserves the contraction
$\iota (X) \omega$.  Therefore for any vector $A$ in the Lie algebra
$\fg$ of $G$ we have $L_{A_M} (\iota (X) \omega) = 0$, where $A_M$ as
before denotes the vector field on $M$ induced by $A$.  Therefore 
$$
0= d \iota (A_M) \iota (X) \omega + \iota (A_M) d (\iota (X) \omega)
= d (\omega (X, A_M)) + \iota (A_M) \omega , 
$$ 
and consequently 
$$
\iota (A_M) \omega = d (\omega (A_M,X)). 
$$
We conclude that the map $\Phi : M \to \fg^*$ defined by 
$$
\langle \Phi (m) , A \rangle = \omega _m (X(m), A_M (m))
$$
is a moment map for the action of $G$ on $(M, \omega)$.
\end{proof}

\begin{exercise}\labell{exercise-beta-mm}
Suppose a Lie group $G$ acts on a manifold $M$ preserving a 1-form $\beta$.  Define the {\bf $\beta$-moment map } $\Psi _\beta : M \to \fg^*$ by
$$
\langle \Psi_\beta (m) , A\rangle = \beta _m (A_M (m))
$$ 
for all $A\in \fg$ and all $m\in M$.  Here as usual $A_M$ denotes
the vector field induced by $A$ on $M$.   

Show that $\Psi_\beta$ is $G$-equivariant, that is, show that 
for any $a\in G$ and any $m\in M$
$$
\Psi_\beta (a\cdot m) = Ad^\dagger (a) \Psi_\beta (m),
$$ where $Ad^\dagger : G \to \text{GL} (\fg^*)$ denotes the coadjoint
representation.  Conclude that the map $\Phi$ defined in
Proposition~\ref{mmap-on-cones} is equivariant.
\end{exercise}

\begin{definition}
A 1-form $\alpha$ on a manifold $B$ is a {\bf  contact form} if the
following two conditions hold:
\begin{enumerate}
\item 
$\alpha _b \not = 0$ for all points $b\in B$.  Hence $\xi: = \ker
\alpha = \{ (b, v)\in TB \mid \alpha _b (v) = 0\}$ is a vector subbundle 
of the tangent bundle $TB$.

\item
$d\alpha |_\xi$ is a symplectic structure on the vector bundle $\xi
\to B$ (i.e. $d\alpha _b |_{\xi_b}$ is nondegenerate).
\end{enumerate}
\end{definition}
\begin{remark}
\begin{enumerate}
\item 
If $\xi \to B$ is a symplectic vector bundle, then the dimension of
its fibers is necessarily even.  Hence if a manifold $B$ has a contact
form, then $B$ is odd-dimensional.
\item
A 1-form $\alpha$ on $2n+1$ dimensional manifold $B$ is contact if and only if
the form $\alpha \wedge (d\alpha)^n$ is never zero, i.e., it is a volume form.
[Prove this].
\end{enumerate}
\end{remark}

\begin{example}
The 1-form $\alpha = dz + x \,dy$ on $\R^3$ is a contact form: $\alpha
\wedge d\alpha = dz \wedge dx \wedge dy$.
\end{example}

\begin{example}\labell{ex2}
Let $B = \R \times \bbT^2$.  Denote the coordinates by $t, \theta_1$
and $\theta_2$ respectively.  The 1-form $\alpha = \cos t \, d\theta_1
+ \sin t \, d\theta_2$ is contact. [Check this.]
\end{example}

\begin{lemma}
Suppose $\alpha$ is a contact form on a manifold $B$. Then for any
positive function $f$ on $B$ the 1-form $f\alpha$ is also contact.
\end{lemma}

\begin{proof}
Note first that since $f$ is nowhere zero, $ker f\alpha = \ker
\alpha$.  Thus to show that $f\alpha$ is contact, it is enough to
check that $d(f\alpha)|_\xi$ is nondegenerate, where $\xi = \ker
\alpha = \ker f\alpha$.  Now 
$d(f\alpha) = df \wedge \alpha + f d\alpha$ and $\alpha |_\xi =0$.
Therefore $d(f\alpha) |_\xi = f d\alpha |_\xi$.  But $f$ is nowhere
zero and $d\alpha |_\xi$ is nondegenerate by assumption.  Thus
$d(f\alpha) |_\xi$ is nondegenerate.
\end{proof}

\begin{definition}
We define the {\bf conformal class} of a 1-form $\alpha$ on a manifold
$B$ to be the set $[\alpha] = \{ e^h \alpha \mid h \in C^\infty (B)\}$,
that is, the set of all 1-forms obtained from $\alpha$ by multiplying
it by a positive function.
\end{definition}
Thus if a 1-form $\alpha$ on a manifold $B$ is contact, then its
conformal class consists of contact forms all defining the same
subbundle $\xi$ of the tangent bundle of $B$.

\begin{definition}\label{def-cont-structure}
A (co-orientable) {\bf contact structure} $\xi$ on a manifold $B$ is a
subbundle of the tangent bundle $TB$ of the form $\xi = \ker \alpha$
for some contact form $\alpha$.

A {\bf co-orientation} of a contact structure $\xi$ is a choice of a
conformal class of contact forms defining the contact structure.
\end{definition}

\begin{remark}
More generally a {\bf contact structure} on a manifold $B$ is a
subbundle $\xi$ of the tangent bundle $TB$ such that for every point
$x\in B$ there is a contact 1-form $\alpha$ defined in a neighborhood
of $x$ with $\ker \alpha = \xi$.  There exist contact structures which
are not co-orientable. For such structures $\xi$ a one-form $\alpha$
with $\ker \alpha = \xi$ exists only {\bf locally}. {\sc In these
notes we will only deal with co-orientable contact structures }.
\end{remark}

\begin{exercise}\label{exercise2.16}
Let $\beta$ be a nowhere zero 1-form on a manifold $B$ and let $\eta =
\ker \beta$.  Let $\eta^\circ \to B$ denote the annihilator of $\eta $ in
$T^*B$: the fiber of $\eta^\circ $ at $b\in B$ is the vector space
$$
\eta^\circ _b = \{ p \in T^*_b B \mid p|_{\eta_b} = 0\}.
$$ 
Show that $\beta$ is a nowhere zero section of real line bundle
$\eta^\circ \to B$.  Show that any other nowhere zero section $\beta'
$ of $\eta ^\circ \to B$ is of the form $\beta' = f \beta$ for some
nowhere zero function $f$ on $B$.

Conclude that a contact structure $\xi$ is co-orientable if and only
if the punctured real line bundle $\xi^\circ \smallsetminus 0$ has two
components (0 of course denotes the zero section).  Show that a choice
of co-orientation of $\xi$ is the same as a choice of a component
$\xi^\circ_+$ of the punctured  bundle $\xi^\circ \smallsetminus 0$.
\end{exercise}

\begin{definition}
Let $(B_1, \xi_1 = \ker \alpha_1 )$ and $(B_2, \xi_2 = \ker \alpha_2)$
be two co-oriented contact manifolds.  A diffeomorphism $\varphi: B_1
\to B_2$ is a {\bf contactomorphism} if the differential $d\varphi$
maps $\xi_1$ to $\xi_2$ preserving the co-orientations.  That is,
$\varphi^* \alpha_2 = f \alpha_1$ for some positive function
$f$.\footnote{ Since $d\varphi (\xi_1) = \xi_2$, the lift
$\tilde{\varphi}: T^*B_1 \to T^*B_2$ of $\varphi$ maps $\xi_1^\circ $
to $\xi_2^\circ$.  ``$d\varphi$ preserves the co-orientation'' means
that $\tilde{\varphi} $ maps $(\xi_1)^\circ_+$ to $(\xi_2)^\circ_+$ (cf.\ Exercise~\ref{exercise2.16}). }
\end{definition}

\begin{definition}
An action of a Lie group $G$ on a manifold $B$ {\bf preserves a contact
structure} $\xi$ and its co-orientation if for every element $a\in G$
the corresponding diffeomorphism $a_B :B \to B$ is a contactomorphism.
We will also say that the action of $G$ on $(B, \xi)$ is a {\bf
contact action}.
\end{definition}
\begin{definition} \labell{def-reeb}
Let $\alpha$ be a contact form on a manifold $B$.  The {\bf Reeb}
vector field $Y_\alpha$ of $\alpha$ is the unique vector field
satisfying $\iota (Y_\alpha) d\alpha = 0$ and $\alpha (Y_\alpha) = 1$.
\end{definition}
\begin{exercise}
Why does the definition of the Reeb vector field makes sense?  
\end{exercise}
\begin{remark}
The Reeb vector field depends strongly on the contact form.  And it is
not just its magnitude: if $\alpha$ is a contact form and $Y_\alpha $
is its Reeb vector field, then there is no reason for $\iota
(Y_\alpha) d(f\alpha) = 0$ where $f$ is a nowhere zero function.
\end{remark}

\begin{exercise}
Compute the Reeb vector field of the contact form $\alpha = dz + x\,
dy$ on $\R^3$.
\end{exercise}

\begin{exercise}
Compute the Reeb vector field of the contact form $\alpha$ of
Example~\ref{ex2}: $\alpha =  \cos t \, d\theta_1 + \sin t \,
d\theta_2$ on $B = \R \times \bbT^2$.
\end{exercise}

Symplectic cones and contact manifolds are intimately related:
\begin{theorem}\labell{thm-cone-cont}
Suppose a compact connected Lie group $G$ acts effectively on a
symplectic cone $(M, \omega, X)$ preserving the symplectic form
$\omega$ and the expanding vector field $X$.  Then $G$ induces an
effective action on the base $B$ of the cone making the projection
$\varpi : M \to B$ $G$-equivariant.  Moreover, the base $B$ has a
natural co-oriented contact structure $\xi$, and the induced action of
$G$ on $B$ preserves a contact form $\alpha$ defining $\xi$.  In
particular the action of $G$ on $(B, \xi)$ is contact.
\end{theorem}
We start the proof of Theorem~\ref{thm-cone-cont} with an observation
that the action of $G$ on $M$ descends to an effective action of $G$
on the base $B$ making the projection $\varpi : M \to B$
$G$-equivariant.  Next we prove:

\begin{proposition}\labell{triv-prop}
Any principal $\R$-bundle $\R \to M \stackrel{\varpi}{\to} B$ is trivial. 
\end{proposition}
\begin{proof}
The proposition is true because the real line is contractible.  Here
is an elementary argument.   Note first that if $s: B
\to M$ is a (local) section of $ \to M \stackrel{\varpi}{\to} B$ and
$f\in C^\infty (B)$ is a function, then $s-f$ makes sense; it is again
a (local) section of $\varpi: M\to B$.  To prove that a principal
bundle is trivial it is enough to construct a global section.  To this
end choose an open cover $\{U_\alpha\}$ of $B$ such that for each
$U_\alpha$ there is a section $s_\alpha : U_\alpha \to M$.  Choose a
partition of unity $\tau_\alpha$ subordinate to the cover
$\{U_\alpha\}$.  Two sections of a principal $\R$-bundle differ by
real-valued function.  Thus by abuse of notation on an intersection
$U_\alpha \cap U_\beta$, $s_\alpha -s_\beta$ is a real-valued
function.  Now define for each index $\alpha$ $$ s_\beta ' = s_\beta -
\sum _{\alpha \not = \beta}
\tau_\alpha (s_\beta - s_\alpha).  $$ Then on an intersection
$U_\alpha \cap U_\beta$
\begin{equation*}
\begin{split}
s_\beta ' - s_\gamma ' 
	& = \left(s_\beta - \sum _{\alpha \not = \beta} 
\tau_\alpha (s_\beta - s_\alpha) \right) -
\left(s_\gamma - \sum _{\alpha \not = \gamma} \tau_\alpha (s_\gamma - s_\alpha)
\right)\\
	&= s_\beta - s_\gamma - \left(\sum _{\alpha \not = \beta, \gamma}  \tau _\alpha  (s_\beta - s_\gamma )\right) + \tau _\beta (s_\gamma - s_\beta) -
\tau_\gamma (s_\beta - s_\gamma) \\
&=s_\beta - s_\gamma - (\sum _\alpha \tau_\alpha ) (s_\beta - s_\gamma) =0 . 
\end{split}
\end{equation*}
Therefore the collection of local sections $\{s_\alpha '\}$ defines a
global section of $\varpi : M \to B$.  Consequently the bundle is trivial.
\end{proof}
Thus  any symplectic cone is of the form $B\times \R$ where $B
= M/\R$ is an odd-dimensional manifold.

\begin{lemma}\labell{another-lemma}
Let $(M, \omega, X)$ be a symplectic cone, let $B$ be its base and let
$\varpi : M \to B$ denote the projection.  Pick a trivialization
$\varphi: B \times \R \to M$.  Then $\varphi^* \omega = d (e^t
\alpha)$ where $t$ is a coordinate on $\R$ and $\alpha $ is a contact
form on $B$.  Conversely, if $\alpha$ is contact form on $B$ then
$(B\times \R, d(e^t \alpha), \frac{\partial}{\partial t})$ is a
symplectic cone.
\end{lemma}

\begin{proof}
By Proposition~\ref{triv-prop} the principal $\R$ bundle $\varpi: M\to
B$ is trivial.  Choose a trivialization $M \simeq B \times \R$.  Under
this identification the vector field $X$ becomes
$\frac{\partial}{\partial t}$.

Since $d\omega =0$ and $L_X \omega = \omega$, $d \iota (X) \omega =
\omega$ (c.f. proof of Proposition~\ref{mmap-on-cones}).  Let 
$\beta = \iota (X) \omega$.  Then $\iota (X) \beta = 0$ and $L_X \beta
= d \iota (X) \beta + \iota (X) d\beta = 0+ \iota (X) \omega = \beta$.
Hence for any point $(b, t) \in B \times \R$
$$
\beta_{(b,t)} = e^t \beta_{(b, 0)}.
$$ 
Since $\iota (X (b, 0))\beta_{(b, 0)} = 0$ it follows that
$\beta_{(b, 0)} = \alpha _b$ for a 1-form $\alpha $ on $B$.  It
remains to show that $\alpha$ is contact.  For this it suffices to
show that $\alpha \wedge (d\alpha)^d $ is nowhere zero, where $d =
\frac{1}{2} \dim M -1$.  Now $\omega = d(e^t \alpha) $ is symplectic.
Hence $\omega^{d+1} $ is nowhere vanishing.  Now $\omega^{d+1} = (e^t
(dt \wedge \alpha + d\alpha))^{d+1} = e^{td} dt \wedge \alpha \wedge
(d\alpha)^d$.  Hence  $\alpha \wedge (d\alpha)^d $ is nowhere vanishing.

Conversely suppose $\alpha$ is a contact 1-form on $B$. Let $\omega =
d(e^t \alpha)$ and let $X = \frac{\partial}{\partial t}$.
Then $L_X \omega = d(\iota (\frac{\partial}{\partial t}) d (e^t
\alpha))= d (\iota (\frac{\partial}{\partial t}) ( e^t dt \wedge \alpha + e^t
d\alpha)) = d (e^t \alpha + 0 ) = \omega$.
It remains to check that $\omega$ is nondegenerate.  For any $(b, t)
\in B\times \R$, the tangent space $T_{(b,t)} (B\times \R)$ decomposes
as $T_{(b,t)} (B\times \R) = \ker \alpha_b \oplus \R Y_\alpha
(b)\oplus \R$ where $Y_\alpha $ is the Reeb vector field of $\alpha$
(cf.\ Definition~\ref{def-reeb}). Since $\alpha$ is contact $d\alpha _b
|_{\ker \alpha_b}$ is nondegenerate.  The restriction $dt \wedge
\alpha_b $ to $\R Y_\alpha (b)\oplus \R$ is nondegenerate as well.
Hence $\omega = e^t ( dt
\wedge \alpha + d\alpha )$ is nondegenerate.  This proves that
$(B\times \R, d(e^t \alpha), \frac{\partial}{\partial t})$ is a
symplectic cone.
\end{proof}

\begin{exercise}\labell{cont-form-exerc}
Let $(M, \omega, X)$ be a symplectic cone, let $B$ be its base and let
$\varpi : M \to B$ denote the $\R$-orbit map. Pick a global section
$s: B \to M$ of $\varpi : M \to B$ and let $\alpha = s^* (\iota
(X)\omega)$.  Show that $\alpha$ is a contact form on $B$.  Show that
it is the same contact form that the proof of
Lemma~\ref{another-lemma} would produce from the trivialization
$\varphi : B \times \R \to M$, $\varphi (b, t) = \rho _t (s(b))$. Here
$\rho_t :M \to M$ denotes the action of $\R$.
\end{exercise}

\begin{remark}
Different choices of trivializations $\varphi$ of $\varpi : M\to B$
give rise to different contact forms on $B$.  However, they all define the
same contact structure $\xi$ on the base $B$.  Intrinsically $\xi$ can
be defined as follows: for a point $b\in B$, 
\begin{equation}\label{eq-xi}
\xi _b = d \varpi_ m (\ker (\iota (X) \omega)_m) 
\quad \text{for any }m\in \varpi \inv (b) 
\end{equation}
It is not hard to check that $\xi $ is well-defined.  First note that
$\R$ acts transitively on the fiber $\varpi \inv (b)$.  Second observe
that for any $\lambda \in \R$ we have $\rho_\lambda ^* (\iota (X)
\omega) = e^\lambda (\iota (X) \omega)$, and hence $d\rho _\lambda
(\ker (\iota (X) \omega)_m) = (\ker (\iota (X) \omega)_{\rho_\lambda
(m)}$.  Here again $\rho_t :M \to M$ denotes the action of $\R$.
\end{remark}
It follows that the action of $G$ on $B$ induced by an action of $G$
on the symplectic cone $\varpi: M\to B$ preserves the contact
structure $\xi$ defined by (\ref{eq-xi}).  Since $G$ is connected and
since the identity map preserves the co-orientation of $\xi$, all the
other elements of $G$ also preserve the co-orientation.

It remains to show that there is a $G$-invariant 1-form $\alpha$ with $\ker
\alpha = \xi$.  By Lemma~\ref{another-lemma} a choice of a trivialization 
of $\varpi : M\to B$ gives us a 1-form $\alpha$ on $B$ with $\ker
\alpha = \xi$, but this form need not be $G$-invariant.  It is only
$G$-invariant if the trivialization is $G$-equivariant.  Therefore we
proceed as follows.

\begin{lemma}\label{lemma2.9}
Suppose a compact Lie group $G$ acts on a manifold $B$ preserving a
(co-oriented) contact structure $\xi = \ker \alpha$ for some 1-form
$\alpha$.  Then there exists a $G$-invariant
1-form $\tilde{\alpha}$ with $\ker
\tilde{\alpha} = \xi$.
\end{lemma}

\begin{proof} For every $a\in G$ the corresponding diffeomorphism 
$a_B : B \to B$ is a contactomorphism.  Hence $(a_B)^* \alpha = f_a
\alpha$ for some positive function $f_a$ depending smoothly on $a$.
Define a new contact form $\tilde{\alpha}$ to be the average of
$\alpha$ over the action of $G$: 
$$
\tilde{\alpha}_b = \int _G \left((a_B)^*\alpha \right)_b \, da  = 
\int _G \left( f_a (b) \alpha _b \right) \, da  = 
\left(\int _G f_a (b) \,da\right) \alpha_b
$$ 
for all $b\in B$. Here $da$ is a bi-invariant measure on $G$
normalized so that $\int _G da = 1$.  Since 
 $f_a > 0$ for all $a\in G$, the integral
$\left(\int _G f_a (b) \,da\right)$ is positive and $\tilde{\alpha}$
is indeed nowhere zero.
\end{proof}
{\it From now on we will always assume that whenever a group actions
preserves a contact structure it also preserves a contact form
defining this structure.}

This concludes the proof of
Theorem~\ref{thm-cone-cont}. It follows from the Theorem that if an
$n$-torus $G$ acts effectively on the punctured cotangent bundle $M =
T^* \bbT^n \smallsetminus 0$ preserving the symplectic form and
commuting with dilations then it acts on the quotient $B = M/\R \simeq
S^* \bbT^n$ preserving the corresponding contact structure.  Note that
$2 \dim G = \dim B + 1$.

\begin{definition}\label{def-ctm}
An effective action of a torus $G$ on a manifold $B$ preserving a contact
structure $ \xi$ is {\bf completely integrable}  if $2\dim G = \dim B +1$.

A {\bf contact toric $G$-manifold} is a co-oriented contact manifold
with a completely integrable action of a torus $G$.
\end{definition}

We are now in position to reduce Theorem~\ref{main} to a statement
about contact toric manifolds.  Consider again the action of
$G=\bbT^n$ on $M = T^* \bbT^n \smallsetminus 0$ preserving the
symplectic form and commuting with dilations.  As was remarked
previously $M$ is a symplectic cone over $B = S^* \bbT^n = \bbT^n
\times S^{n-1}$.  By Theorem~\ref{thm-cone-cont} the action of $G$ on
$M$ induces an effective action on $B$.  Moreover $B$ has a
$G$-invariant contact structure $\xi$ making $(B, \xi)$ into a compact
connected contact toric $G$-manifold.  Clearly if the action of $G$ on
$B$ is free then the original action of $G$ on $M$ was free as well.
Therefore a proof of Theorem~\ref{main} reduces to
\begin{theorem}\label{main-thm'}
Let $(B, \xi = \ker \alpha)$ be a compact
connected contact toric $G$-manifold and suppose the action of $G$ is
not free.  Then $B$ is not homotopy equivalent to $\bbT^n \times
S^{n-1}$, $n = \dim G$.
\end{theorem}

The rest of the notes will be occupied with a proof of
Theorem~\ref{main-thm'}.  The proof uses heavily contact moment maps
which are discussed in the next section.  We end this section with a
partial converse to Theorem~\ref{thm-cone-cont}.

We have seen that given a contact form $\alpha$ on a manifold $B$ the
form $d( e^t \alpha)$ on $B\times \R$ is symplectic.  The pair
$(B\times \R,d( e^t \alpha) )$ is called the {\bf
symplectization}\footnote{Sometimes $(B\times \R,d( e^t \alpha) )$ is
called the {\bf symplectification} of $(B, \alpha)$.}  of $(B,
\alpha)$.  It is clearly a symplectic cone.

Different contact forms on $B$ give rise to different symplectic forms
on $B\times \R$.  However there is a symplectic cone that depends only
on the contact structure $\xi = \ker \alpha$ (and its co-orientation)
and not on a particular choice of a contact form: Let $\xi_+^\circ$
denote the component of $\xi^\circ \smallsetminus 0$ giving $\xi $ its
co-orientation (cf.\ Exercise~\ref{exercise2.16}).  It is not hard to
check that $\xi^\circ_+$ is a symplectic submanifold of the cotangent
bundle $T^*B$ with its standard symplectic structure.  The action of
$\R$ on $T^*B \smallsetminus 0$ by dilations preserves $\xi^\circ_+$
and makes it into a symplectic cone.  A section $\beta$ of
$\xi_+^\circ \to B$ is a contact form $\beta$ on $B$ with $\xi = \ker
\beta$.  Moreover the trivialization $\varphi_\beta :B \times \R \to
\xi^\circ_+$, $\varphi_\beta (b, t) = e^t \beta_b$, that $\beta$
defines pulls back the symplectic form on $\xi^\circ_+$ to $d(e^t
\beta)$.

Finally given a symplectic cone $(M, \omega, X)$ with the base $B$, the
orbit map $\varpi : M \to B$ and the induced contact structure $\xi$ on
$B$, there is an $\R$ equivariant symplectomorphism $\varphi: M \to
\xi^\circ_+$ defined as follows: for  a point $m\in M$ let $\varphi (m)$ be  
the covector in $T^*_{\varpi(m)} B$ such that
$$
\varphi (m) (v) = (\iota (X)\omega)_m ( ds_b (v))
$$ 
for all $v\in T_{\varpi (m)} B$ and a  section $s$ of $\varpi
:M \to B$.

The discussion above can be summarized as:

\begin{lemma}
Let $(B, \xi = \ker \alpha)$ be a contact manifold, let $\xi^\circ_+$
be a component of $\xi^\circ \smallsetminus 0$, the annihilator of
$\xi$ in $T^*B$ minus the zero section.  The principal $\R$ bundle
$\xi^\circ_+ \to B$ is a symplectic cone.

If $(M, \omega, X)$ is a symplectic cone with the base $B$ and $\xi$
is the induced contact structure on $B$, then $\xi^\circ_+$ is
isomorphic to $(M, \omega, X)$ as a symplectic cone.
\end{lemma}

\chapter{Group actions and moment maps on contact manifolds}

Moment maps exist in the category of contact group actions.  In fact
moment maps exist for all contact actions.  This is because a contact
form defines a bijection between contact vector fields and smooth
functions.

\begin{definition}
A vector field $X$ on a contact manifold $(B, \xi = \ker \alpha)$ is
{\bf contact} if its flow $\varphi_t$ consists of contactomorphisms.
In particular $d \varphi_t (\xi) \subset \xi$. Hence for any section
$v$ of the bundle $\xi \to B$, the Lie bracket $[X, v]$ is again a
section of $\xi \to B$.
\end{definition}
Thus for a contact action of a Lie group $G$ on $(B, \xi)$ the vector
fields induced by elements of the Lie algebra $\fg$ of $G$ are contact.

\begin{exercise}\label{exercise3.2}
Prove that a Reeb vector field is contact.  More generally prove that
a vector field $X$ on a contact manifold $(B, \xi = \ker \alpha)$ is
contact if and only if $L_X \alpha = h \alpha$ for some function $h$
($h$ can have zeros).
\end{exercise}

A choice of a contact form on a contact manifold $(B, \xi)$ identifies
contact vector fields with smooth functions.

\begin{proposition}
Let $(B, \xi = \ker \alpha)$ be a contact manifold.  The linear map from
contact vector fields to smooth functions given by
$X \mapsto f^X := \alpha (X)$ is one-to-one and onto.
\end{proposition}

\begin{proof}
Note that the Reeb vector field $Y_\alpha$ corresponds to the function
1.  For any vector field $X$ on $B$ the vector field $X - \alpha (X)
Y_\alpha$ is in the kernel of $\alpha$, which is the contact
distribution $\xi$.  Since $d\alpha |_\xi$ is non-degenerate, $X -
\alpha (X) Y_\alpha$ is uniquely determined by $\iota (X - \alpha (X)
Y_\alpha) (d\alpha |_\xi)$.

For any section $v$ of $\xi \to B$ and any vector field $X$ we have 
$$ 
0= L_X 0 = L_X (\alpha (v)) = (L_X \alpha) (v) + \alpha ([X, v]).  
$$ 
Assume now that $X$ is contact.  Then $\alpha ([X, v]) = 0$.  Therefore,
by Cartan's formula
$$
0 = (d\iota (X)\alpha + \iota (X) d\alpha ) (v).
$$
Hence for any section $v$ of $\xi$, $d\alpha (X, v) = - d(\alpha (X))
(v) = -d f^X (v)$.  And, of course, $d\alpha (X, v) = d \alpha (X -
f^X Y_\alpha, v)$ for all $v$. We conclude that 
\begin{equation}\label{eq*}
\iota (X - \alpha (X) Y_\alpha) (d\alpha |_\xi) = -df^X |_\xi
\end{equation}
for any contact vector field $X$.  Thus if $X$ is contact, its
component in the direction of the Reeb vector field is $f^X Y_\alpha$
and its component in the direction of the contact distribution is
uniquely determined by (\ref{eq*}).

Conversely, given a function $f$ on $B$ there is a unique section
$X'_f$ of $\xi$ such that
\begin{equation}\labell{eq2}
\iota (X'_f) d\alpha |_\xi = - df |_\xi.
\end{equation}
The vector field $X_f := X'_f + f Y_\alpha$ is a contact vector
field with $\alpha (X_f) = f$. 
\end{proof}

\begin{exercise}
Given a function $f$ on a manifold $B$ with contact form $\alpha$
check that the vector field $X_f := X'_f + f Y_\alpha$, where $X'_f$
is defined by (\ref{eq2}), is contact.  That is, show that $L_{X_f}
\alpha = h \alpha$ for some function $h$ (cf.\ Exercise~\ref{exercise3.2}).
\end{exercise}

Suppose now that a Lie group $G$ acts on a manifold $B$ preserving a
contact 1-form $\alpha$.  Then for any vector $A$ in the Lie algebra
$\fg$ of $G$, the induced vector field $A_B$\footnote{ Recall that the
 vector field $A_B$ is
defined by $A_B (b) = \left. \frac{d}{dt} \right|_{t= 0} \exp (tA)
\cdot b$.}  satisfies $L_{A_B} \alpha = 0$ and, in particular, is contact.
Since the contact form $\alpha$ defines a 1-1 correspondence between
contact vector fields and functions, it makes sense to define the {\bf
$\alpha$-moment map} $\Psi_\alpha : B \to \fg^*$ by
\begin{equation}
\langle \Psi _\alpha (b) , A\rangle = \alpha _b (A_B (b))
\end{equation}
for all $b\in B$ and all $A\in \fg$.  Note that by
Exercise~\ref{exercise-beta-mm} the $\alpha$-moment map $\Psi_\alpha$
is $G$-equivariant.

The $\alpha$-moment map, as the name suggests, depends rather
strongly on $\alpha$: if $f$ is any positive $G$-invariant function on
$B$, then $f\alpha$ is another $G$ invariant contact form and clearly
$$
\Psi_{f\alpha} = f\Psi_\alpha.
$$ In particular, unlike in the symplectic case, the image of a moment
map is {\bf not} an invariant of the action and of the contact
structure. However the {\bf moment cone} $C(\Psi _\alpha)$ defined by
$$ 
C(\Psi _\alpha) = \{ t\eta \in \fg^* \mid \eta \in \Psi_\alpha (B),
\, t\in [0, \infty) \}
$$
is an invariant of the action of $G$ on $(B, \xi = \ker \alpha)$.

\begin{remark}
Suppose a Lie group $G$ acts on a manifold $B$ preserving a contact
form $\alpha$.  The action of $G$ lifts to an action on the cotangent
bundle $T^*B$ which preserves the annihilator $\xi^\circ $ of the
contact structure $\xi = \ker \alpha$.  Moreover, the action preserves
the component $\xi^\circ_+$ of $\xi \smallsetminus 0$ which contains
the image of $\alpha :B \to T^*B$.

The action of $G$ on $T^*B$ is Hamiltonian with a natural moment map
$\Phi: T^*B \to \fg^*$ given by
$$
\langle \Phi (b, p), A \rangle = \langle p, A_B (b)\rangle
$$ for all $b\in B$, $p\in T^*_b B$, $A\in \fg$.  Since the
submanifold $\xi^\circ_+$ is a $G$-invariant symplectic submanifold of
$T^*B$ the restriction $\Psi:= \Phi |_{\xi^\circ_+}$ is a moment map for the
action of $G$ on $\xi^\circ_+$.  It is not hard to check that
$$
\Psi_\alpha = \Psi \circ \alpha .
$$ 
For this and other reasons it makes sense to think of $\Psi: \xi^\circ_+
\to \fg^*$ as the moment map for the action of $G$ on 
$(B, \xi)$.\footnote{Note that since $\Phi$ is $G$-equivariant and
$\alpha : B \to T^*B$ is $G$-equivariant, $\Psi_\alpha$ is
$G$-equivariant as well. This gives us an alternative proof that the
$\alpha$-moment map is equivariant.}  Note also that $C(\Psi_\alpha) =
\Psi (\xi^\circ _+)
\cup \{0\}$.  We will thus denote the moment cone $C(\Psi_\alpha)$ by
$C(\Psi)$.
\end{remark}

Consequently and by analogy with symplectic toric manifolds we will
think of {\bf contact toric $G$-manifolds} as triples $(B, \xi = \ker
\alpha,\Psi: \xi^\circ_+ \to \fg^*)$.

\chapter{Contact toric manifolds}

The rest of the lecture notes will be devoted to a proof of
Theorem~\ref{main-thm'}.  Right from the beginning the proof will
bifurcate into two cases: the contact manifold $B$ is 3-dimensional
and $\dim B > 3$.  If $\dim B = 3$ we will argue directly using slices
that the orbit space $B/G$ is homeomorphic to a closed interval $[0,
1]$ and then use this to compute the integral cohomology of $B$.  This
will show that $B$ cannot be homeomorphic to $S^*\bbT^2 = \bbT^3$.

We will then consider the case where $\dim B > 3$.  In this case we have a
connectedness and convexity theorem of Banyaga and Molino (see \cite{BM1, BM2};
for a different proof see \cite{L-conv}):

\begin{theorem}\label{thm-BM}
Let $(B, \xi = \ker \alpha, \Psi: \xi^\circ_+ \to \fg^*)$ be a compact
connected contact toric $G$-manifold.  Suppose $\dim B > 3$.  Then the
fibers of the moment map $\Psi$ are $G$-orbits (and in
particular are connected) and the moment cone $C(\Psi)$ is a
convex polyhedral cone in $\fg^*$.  Moreover $C(\Psi) \not =
\fg^*$ iff the action of $G$ is not free.
\end{theorem}

Our proof will then bifurcate again. We will consider separately the case where
the moment cone contains a linear subspace of dimension $k$, $0 < k <
\dim G$ and where no such subspace exists (i.e., the moment cone is {\bf
proper}).  

In the first case we will use a uniqueness theorem of Boucetta and
Molino \cite{BoM}\footnote{The result was rediscovered a few years
later by Lerman, Tolman and Woodward (Lemma~7.2 and Proposition~7.3 in
\cite{LT}).} for symplectic toric manifolds to argue that the
symplectization $B\times \R$ of $B$ is diffeomorphic to the manifold
$N = \bbT^k \times (\R^k \times \C^l \smallsetminus \{(0,0)\})$ where
$2l + 2 k = \dim B + 1$ and $k, l> 0$.  It is easy to see that $N$
cannot be homotopy equivalent $T^* \bbT^n \smallsetminus 0$, $n = k
+l$.

In the latter case we will argue following Boyer and Galicki \cite{BG}
that there is a locally free\footnote{An action of a Lie group $G$ on
a manifold $Z$ is locally free if all the isotropy groups are zero
dimensional.} $S^1$ action on $B$ such that the quotient $B/S^1$ is a
(compact connected) symplectic toric orbifold.  Since the action of
$S^1$ is locally free there is a long exact sequence of rational
cohomology groups (the Gysin sequence) tying together the cohomology
of $B$ and of $B/S^1$.On the other hand Morse theory on the orbifold
$B/S^1$ shows that all odd-dimensional rational cohomology of $B/S^1$
vanishes.  Together these two facts will imply that $\dim _\Q H^1 (B,
\Q) \leq 1$ .  Thus $B$ cannot be $S^* \bbT^n = \bbT^n \times
S^{n-1}$, $n = \frac{1}{2} (\dim B + 1 ) > 2$.  This completes the
preview of our proof of Theorem~\ref{main-thm'}.

\section{Homogeneous vector bundles and slices}

Suppose that $G$ is a Lie group, $H \subset G$ a closed subgroup and
suppose we have a representation of $H$ on a vector space $W$.  Then
$H$ acts on the product $G\times W$ by $h\cdot (g, w) = (gh\inv,
g\cdot w)$ for $h\in H$, $(g, w) \in G\times W$.  The quotient
$G\times _H W := (G \times W)/H$ is a vector bundle over $G/H$ with
typical fiber $W$.  We denote the image of $(g, w) \in G\times W$ in
$G\times_H W$ by $[g, w]$.  Note that the action of $G$ on $ G\times
W$ given by $a\cdot (g, w) = (ag, w)$ commutes with the action of $H$
and hence descends to an action $G$ on $G\times _H W$: $a\cdot [g, w]
= [ag, w]$.  The projection $G\times _H W \to G/H$ is $G$-equivariant
and the action of $G$ on the base $G/H$ is transitive.  This makes
$G\times _H W \to G/H$ into a homogeneous vector bundle.

Conversely if $\pi : E \to G/H$ is a vector bundle with an action of
$G$ by vector bundle maps making $\pi$ equivariant, then $E$ is
isomorphic to $G\times _H W$ where $W$ is the fiber of $E$ above the
identity coset $1H$.\footnote{Why is there a representation of $H$ on
$W$?}  Indeed the map $$ G \times W \to E, \quad (g, w) \mapsto g
\cdot w $$ is onto and is constant along the orbits of $H$. It
descends to a vector bundle isomorphism 
$$ 
G \times_H W \to E, \quad [g, w] \mapsto g \cdot w .  
$$

Next suppose a compact Lie group $G$ acts on a manifold $M$.  Consider
an orbit $G\cdot x \subset M$.  The group $G$ acts on the normal
bundle of the orbit $\nu (G\cdot x) = TM|_{G\cdot x}/T(G\cdot x)$
making the projection $\pi : \nu (G\cdot x) \to (G\cdot x)$
equivariant.  Thus $ \nu (G\cdot x) = G \times_H W $ where $W = T_x M/
T_x (G\cdot x)$.  

If we choose a $G$-invariant Riemannian metric on $M$ \footnote{Choose
any metric on $M$ and then average it over the group $G$ (cf.\ proof of
Lemma~\ref{lemma2.9}).} we can identify $\nu (G\cdot x)$ with the
perpendicular of $T(G\cdot x)$ in $TM|_{G\cdot x}$.  Furthermore, the
Riemannian exponential map $\exp :
\nu (G\cdot x) \to M$ is $G$-equivariant.  Hence, by the Tubular
Neighborhood theorem we get:

\begin{lemma} \labell{slice-lemma}
Let $G$ be a compact connected Lie group acting on a manifold $M$ and
let $x\in M$ be a point.  A neighborhood of $G\cdot x$ in $M$ is
$G$-equivariantly diffeomorphic to a neighborhood of the zero section
of the homogeneous vector bundle $G\times _{G_x} W$ where $G_x$
denotes the isotropy group of $x$ and  $W = T_x M/ T_x (G\cdot x)$.
\end{lemma}

\begin{definition}
Let $G$ be a Lie group acting on a manifold $M$; let $x\in M$ be a
point.  Denote the isotropy group of $x$ by $G_x$.  An embedded
submanifold $S\subset M$ is a {\bf slice through $x$ for the action of
$G$ on $M$} if $x\in S$, $S$ is $G_x $-invariant and if the map $G
\times S \to M$ given by $(g, s) \mapsto g\cdot s$ descends to an open
embedding $G\times _{G_x} S \to M$, $[g, s] \mapsto g\cdot s$.
\end{definition}

Thus by Lemma~\ref{slice-lemma} slices exist for actions of compact
Lie groups: we may choose as a slice at $x$ the image of a small
$G_x$-invariant neighborhood of 0 in $W = T_x M/T_x (G\cdot x)$ under
the exponential map.

Here is a typical application of the existence of slices. Locally near
$x$ the quotient $M/G$ is homeomorphic to the quotient $(G\times _{G_x}
W)/ G = W/G_x$. Hence quotients of manifolds by actions of compact Lie
groups are modeled on quotients of vector spaces by linear actions of
compact Lie groups.  The linear action of $G_x$ on $W = T_x M /T_x
(G\cdot x)$ is called the {\bf slice representation} at $x$.

Here is another application of the above construction:
\begin{lemma}\labell{lemma-faithful}
Suppose $G$ is a compact abelian group acting effectively on a connected 
manifold $M$.  Then every slice representation is {\bf faithful}, i.e., no
slice representation has a  kernel.
\end{lemma}

\begin{proof} Suppose the slice representation of $H = G_x$ on $W =
T_x M/ T_x (G\cdot x)$ is not faithful at a point $x \in M$.  By
Lemma~\ref{slice-lemma} it is no loss of generality to assume that a
neighborhood of $G\cdot x$ in $M$ is the homogeneous vector bundle
$G\times _H W \to G/H = G\cdot x$ and that the point $x$ is $[1, 0]
\in G\times _H W$.  If there is an element $a\in H$ such that $a\not =
1$ and yet $a\cdot w = w$ for all $w\in W$, then $a\cdot [g, w] = [ag,
w] = [ga, w] = [g, a\cdot w] = [g, w]$ for all $[g, w]
\in G\times _H W$ ($ag= ga$ since $G$ is abelian).  
Thus $a \in H$ fixes an open neighborhood of $x$.  
Since the set fixed by $a$ is closed and since $M$ is connected it
follows that $a$ fixes all of $M$.  This contradicts the assumption
that the action of $G$ on $M$ is effective.
\end{proof}

\begin{remark}\label{rmrk4.7} 
The compactness of the Lie group $G$ is not necessary for the
existence of slices.  According to Palais \cite{P} it is only
necessary that its action on a manifold $M$ be proper (see
Definition~\ref{def-proper} above).

Thus if an action of a Lie group $G$ on a manifold $M$ is free and proper,
then the existence of slices tell us that a neighborhood of every
orbit is equivariantly diffeomorphic to a product of $G$ with some
manifold $S$.  It is not hard to deduce from this that the orbit space
$M/G$ is a manifold and that the orbit map $M \to M/G$ makes $M$ into
a principal $G$-bundle.
\end{remark}

We now recall a few properties of tori, which for us are compact
abelian Lie groups.  If $G$ is a torus, then the (Lie group)
exponential map $\exp :\fg \to G$ is a covering map.  The kernel $\Z_G
$ of $\exp$ is called the {\bf integral lattice}.  Clearly $G =
\fg/\Z_G$.  The group $\Z_G$ is isomorphic to the fundamental group of
$G$. Also it has the property that for any $X\in\Z_G$ the
corresponding 1-parameter subgroup $\{
\exp tX \mid t\in \R\}$ is a circle. The dual lattice 
$\Z_G^* = \text{Hom}_\Z (\Z_G, \Z) \cong 
\{\ell \in \fg^* \mid \ell (\Z_G) \subset \Z\} $ is the weight lattice.  It
parameterizes 1-dimension complex representations of $G$, or,
equivalently, group homomorphisms ({\bf characters}) $\chi : G\to S^1$:
Given $\nu \in \Z_G^*$ the corresponding character $\chi_\nu : G =
\fg/\Z_G \to S^1$ is defined by $\chi_\nu (\exp X) = e^{2\pi i \nu
(X)}$ for all $X \in \fg$.  Given a character $\chi: G \to S^1$, its
differential $d\chi_1$ at 1 is a weight.

Recall that a complex representation of a compact abelian
group is a direct sum of one-dimensional complex representations.  Thus
a complex representation of a torus is completely characterized by a
finite set of weights.  The same is true for {\em symplectic }
representations of tori --- one defines weights with respect to some
complex structure compatible with the symplectic form.  The weights do
not depend on the choice of the complex structure; they only depend on
the symplectic form.

The Lemma below is a key representation-theoretic fact in the
classification of symplectic and contact toric manifolds.
\begin{lemma}\labell{Del-lemma}
Suppose $\rho : H \to \text{Sp} (V, \omega)$ is a symplectic
representation of a compact abelian Lie group $H$ on a symplectic
vector space $(V, \omega)$. Suppose that $\rho$ is faithful, i.e.,
$\ker \rho = \{1\}$.  

Then $\dim H \leq \frac{1}{2} \dim V$.  If $H = \frac{1}{2} \dim V$,
then $H$ is connected.  Moreover, the set of weights for the
representation of $H$ on $(V, \omega)$ is a basis of the weight
lattice $\Z^*_H$ of $H$.
\end{lemma}

\begin{proof}
Since $H$ is compact there exists an $H$-invariant complex structure
$J$ on $V$ compatible with $\omega$ (i.e., $\omega (J \cdot , J\cdot)
= \omega (\cdot , \cdot)$ and for any $v\not = 0$ we have $\omega (Jv,
v) > 0$).

The choice of $J$ identifies $(V, \omega)$ with $\C^n$ with the
standard symplectic structure $\sqrt{-1} \sum d z_j \wedge
d\bar{z}_j$, $n = \frac{1}{2} \dim _{\R} V$.  This, in turn, gives us
a representation of $H$ on $\C^n$ by unitary matrices.  Thus we may
assume that $\rho$ is an injective group homomorphism $\rho :H \to
U(n)$.

The connected component $H^\circ$ of $H$ is a torus. Therefore $\rho
(H^\circ)$ is contained in a maximal torus of $U(n)$ which is the
$n$-torus.  Since $\rho$ has no kernel we have $\dim H = \dim H^\circ
= \dim \rho (H^\circ) \leq n=  \frac{1}{2} \dim V$.

Now suppose $\dim H = n$.  Since all maximal tori in $U(n)$ are
conjugate, we may assume that $\rho(H^\circ)$ is the standard maximal
torus, that is $\rho(H^\circ)$ is the set of all diagonal unitary
matrices. Since the only unitary matrices which commute with all the
diagonal matrices are the diagonal matrices, we see that we must have
$\rho(H) = \rho(H^\circ)$.  Consequently since $\rho$ is faithful, $H=
H^\circ$, i.e., $H$ is a torus.  

Finally the set of weights of the maximal torus $\bbT^n$ in $U(n)$ for
its representation on $\C^n$ is a lattice basis of weight lattice
$\Z_{\bbT^n}^*$.  Hence the set of weights for the representation of
$H$ on $(V, \omega)$ is a basis of the weight lattice $\Z^*_H$ of $H$.
\end{proof}

\begin{lemma}\labell{lemma-noG}
Let $(B, \xi = \ker \alpha)$ be a contact toric $G$-manifold.  Then
\begin{enumerate}
\item No $G$-orbit is tangent to the contact structure $\xi$. In particular 
there are no fixed points.  Hence the $\alpha$-moment map $\Psi_\alpha $ does
not vanish at any point for any $G$-invariant contact form
$\alpha$. Equivalently $\Psi (\xi^\circ_+)$ does not contain the
origin in $\fg^*$.

\item All isotropy groups are connected.
\end{enumerate}
\end{lemma}

\begin{proof}
Let $b$ be a point in $B$ and let $H$ denote its isotropy group.  The
group $H$ acts on the tangent space $T_b B$.  Since the contact form
$\alpha$ is $G$-invariant, its kernel at $b$, the hyperplane $\xi_b$,
is an $H$-invariant subspace of $T_b B$.  Since the Reeb vector field
$Y_\alpha$ of $\alpha$ is unique, the vector $Y_\alpha (b)$ is fixed
by $H$.  Thus we have an $H$ -equivariant splitting 
$$
 T_bB =\R Y_\alpha (b) \oplus \xi _b .
$$

Let $V = T_b (G\cdot b) \cap \xi _b$.  It is an $H$-invariant subspace
of $\xi _b$.  Note that since $\dim \xi_b = \dim B -1$ we either have
$V = T_b (G\cdot b) $ or $\dim V = \dim G\cdot b -1$.  We will argue
that the former case cannot occur. But first we argue that $V$ is an
isotropic subspace of the symplectic vector space $(\xi_b, \omega =
d\alpha_b |_{\xi_b})$.

Let $x, z \in V$ be two vectors. There exist vectors $X, Z \in \fg$
such that $X_B (b) = x$ and $Z_B (b) = z$.  Then $\omega (x, z) =
d\alpha _b (X_B (b), Z_B (b))$.  Since $G$ is abelian $[ Z_B, X_B] =
0$.  The function $\alpha (X_B)$ is $G$-invariant, hence $Z_B (\alpha
(X_B))= 0$.  Therefore
$$
\omega(x,y) =d\alpha (X_B, Z_B) = X_B (\alpha (Z_B)) - Z_B (\alpha (X_B)) - 
d\alpha ([X_B, Z_B]) = 0 - 0 + 0.
$$
This proves that $V$ is isotropic.  

Since $H$ is compact, there is an $H$-invariant complex structure $J$
on $\xi_b$ compatible with $\omega$.  Since $V$ is isotropic, $V\cap
JV = 0$ and $V+ JV = V\oplus JV$ is a symplectic subspace of $(\xi_b ,
\omega)$.  It is $H$-invariant.  In fact, since $G$ is abelian, the
action of $H$ on $T_b (G\cdot b)$ is trivial.  Hence the action of $H$
on $V\oplus JV$ is trivial as well.
 Let $W$ denote the symplectic perpendicular to $V\oplus JV$.  We get a
 symplectic representation of $H$ on $W$.

We now argue that if $T_b (G\cdot b) \subset \xi_b$ then $\dim W$ is
less than $2\dim H$ and that the representation of $H$ on $W$ must be
faithful.  This by Lemma~\ref{Del-lemma} would give us a
contradiction.  We would then conclude that $\dim W = 2 \dim H$, which
by Lemma~\ref{Del-lemma} implies that $H$ is connected.

Suppose now that $T_b (G\cdot b) \subset \xi_b$. Then
 $\dim V = \dim G\cdot b$. Since $B$ is toric, $\dim B = 2 \dim G -1$.  Since 
\begin{equation}\label{H-rep}
T_b B = \R Y_\alpha (b) \oplus V \oplus JV \oplus W
\end{equation}
$\dim W = (2 \dim G - 1 ) - 1 - 2 \dim G\cdot b = 2 \dim G - 2 (\dim G
- \dim H) = 2 \dim H -2$.  Since (\ref{H-rep}) is a splitting as
$H$-representations, the slice representation of $H$ at $b$ is $\R
Y_\alpha (b) \oplus JV \oplus W$.  As we observed earlier the action
of $H$ on $\R Y_\alpha (b) \oplus JV $ is trivial. Since the action of
$G$ on $B$ is effective, the representation of $H$ on $W$ must be
faithful.  Contradiction.

Therefore $T_b (G\cdot b) \not \subset \xi _b$ and so $\dim V = \dim
G\cdot b -1$.  In this case the dimension count gives us exactly that
$\dim W = 2 \dim H$.
\end{proof}

\begin{lemma}\labell{lemma-4.9}
Let $H \subset \bbT^2$ be a closed subgroup isomorphic to $S^1$.  Then
there is another closed subgroup $K \subset \bbT^2$ isomorphic to $S^1$
such that $\bbT^2 = K \times H$.
\end{lemma}

\begin{proof}
Since $H$ is isomorphic to $S^1$ it is of the form $\{\exp t\nu \mid
t\in \R\}$ for some  vector $\nu = (n_1, m_1) \in \Z^2 = \ker \{\exp: \R^2
\to \bbT^2 \}$.  We may assume that $n_1$ and $m_1$ are relatively prime.
Therefore there exist integers $n_2, m_2$ such that $n_1 m_2 - m_1 n_2
= 1$.  Hence the vectors $(n_1, m_1)$ and $(n_2, m_2)$ form a basis of
$\Z^2$. Take $K = \{ \exp t(n_2, m_2) \mid t \in \R\}$.
\end{proof}

\begin{remark} \label{splitting-tori}
More generally if $G$ is a torus and $H \subset G$ is a closed
connected subgroup there is another closed connected subgroup $K
\subset G$ so that $G = K \times H$.  This is a bit harder
to prove than the Lemma above.
\end{remark}

\chapter{\protect Proof of Theorem~\ref{main-thm'} part I: 
the 3-dimensional case} 

In this section we prove Theorem~\ref{main-thm'} in the case that 
$\dim B = 3$:
\begin{lemma}\label{lemma-lens}
Let $B$ be a compact connected contact toric $G =\bbT^2$ manifold (in
particular $\dim B = 3$).  Suppose the action of $G$ is not free.
Then there exist two closed subgroups $K_1, K_2 \subset G$ isomorphic
to $S^1$ so that $B$ is homeomorphic to $([0, 1] \times G)/\sim$ where
$(0, g) \sim (0, ag)$ for all $g\in G$, $a \in K_1$ and $(1, g) \sim
(1, ag) $ for all $g\in G$ and $a\in K_2$.  In other words $B$ is
obtained from the manifold with boundary $[0,1] \times G$ by
collapsing circles in the two components of the boundary by the
respective actions of two circle subgroups. (It may happen that
$K_1 = K_2$).
\end{lemma}

\begin{exercise}
Consider the standard sphere $S^3 = \{ (z_1, z_2) \in C^2 \mid |z_1|^2
+ |z_2|^2 = 1\}$.  The torus $G = \{ (\lambda _1, \lambda_2) \in C^2
\mid |\lambda_1|^2 = 1, \quad |\lambda_2|^2 = 1\}$ acts on $S^3$ by
$(\lambda _1, \lambda _2) \cdot (z_1, z_2 ) = (\lambda _1 z_1, \lambda
_2 z_2)$.  Show that $S^3$ is of the form $([0, 1]\times G)/\sim$ for
two circle subgroups $K_1$, $K_2$ of $G$ where $\sim$ is the
equivalence relation in Lemma~\ref{lemma-lens}.  What are the
subgroups $K_1$, $K_2$?
\end{exercise}

\begin{proof}[Proof of Lemma~\ref{lemma-lens}]
By Lemma~\ref{lemma-noG} all isotropy groups for the action of $G$ on
$B$ are connected and no isotropy group is all of $G$.  Therefore,
since $\dim G = 2$ the possible isotropy groups are trivial or
circles. And points with circle isotropy groups must exist since the
action is not free.

If the isotropy group of a point $b\in B$ is trivial then by
Lemma~\ref{slice-lemma} and the dimension count a $G$-invariant
neighborhood $U$ of $G$ is equivariantly diffeomorphic to $G\times I$
where $I$ is an open interval and $G$ acts on $G\times I$ by $g \cdot
(a, t) = (ga, t)$. Hence $U/G = I$ and there is a map $s: U/G \to B$
so that $\pi \circ s = id$ where $\pi: B \to B/G$ is the orbit map,
i.e., $\pi$ has a local section.

Now consider a point $b\in B$ with the isotropy group $G_b$ isomorphic
to $S^1$.\footnote{By Lemma~\ref{lemma-noG} there are no more
possibilities for $G_b$.}  By Lemma~\ref{slice-lemma} a neighborhood
of $b$ in $B$ is $G$-equivariantly diffeomorphic to a neighborhood of
the zero section in $G\times _{G_b} W$.  By Lemma~\ref{lemma-4.9}
there is a circle $L\subset G$ such that $G= G_b \times L$.  Note that
$G\cdot b \simeq L$ and that $\dim W = 2$ and that by
Lemma~\ref{lemma-faithful} the representation of $G_b$ on $W$ is
faithful.  There is only one faithful real 2-dimensional
representation of $S^1$ --- it is the representation of $S^1$ on
$\R^2$ as $SO(2)$ or, equivalently, the representation of $S^1$ on
$\C$ as $U(1)$.  Since $G=L\times G_b$, $G\times _{G_b} W = L \times
W$, where $L \times G_b$ acts on $L\times W$ by $(\lambda, \mu) \cdot
(g, w) = (\lambda g, \mu \cdot w)$.  We conclude that there is a
$G$-invariant neighborhood $U$ of $G\cdot b$ in $M$, an isomorphism
$\varphi : G\to S^1 \times S^1$ and a diffeomorphism $\psi : U \to S^1
\times D^2$ such that $\psi (g\cdot x ) = \varphi (g) \cdot \psi (x)$
for all $g\in G$, $x\in U$.  Here $(\lambda, \mu)\in S^1 \times S^1
\subset \C \times \C$ acts on $S^1
\times D^2 \subset \C \times \C$ by
$$ 
(\lambda, \mu) \cdot (z_1, z_2) = (\lambda z_1, \mu z_2).  
$$ 
Note
that $S^1 \times D^2 = (\bbT^2 \times [0, 1))/\sim$ where $(\lambda
_1, \lambda _2, 0) \sim (\lambda _1, \mu \lambda _2, 0)$ for all
$(\lambda_1, \lambda _2) \in \bbT^2$ and $\mu \in S^1$.
We conclude that
\begin{enumerate}
\item $U/G \simeq [0, 1)$;

\item $U \simeq (G\times [0,1))/\sim$ where $(g,0) \sim (ag, 0)$ for 
all $g\in G$, $a\in G_b$;

\item there is a local section $s: U/G \to B$ of the orbit map 
$\pi: B \to B/G$;

\item the set of points in $U$ with non-trivial isotropy groups is 
$G\cdot b= \pi\inv (0)$, where again $\pi :U \to [0,1)$
 denotes the orbit map.
\end{enumerate}

It follows that the orbit space $B/G$ is locally homeomorphic to
either an open interval or to a half-open interval $[0,1)$.  Hence
$B/G$ is a topological 1-dimensional manifold with boundary.  Since
$B/G$ is compact and connected we may identify it with $[0, 1]$.  Note
that the set of points with non-trivial isotropy groups is $\pi\inv(
\{0, 1\})$, where by abuse of notation $\pi: B \to [0, 1]$ denotes the
orbit map.  More specifically $\pi\inv (0) = G/K_1$, $\pi\inv (1) =
G/K_2$ for some circle subgroups $K_1, K_2 \subset G$.

We now argue that $\pi:B \to [0,1]$ has a global section $s:[0, 1] \to
B$, so that $\pi \circ s (t) = t $ for all $t\in [0, 1]$.  We have
seen that sections of $\pi$ exist locally: for every $t\in [0, 1]$
there is an interval $I \subset [0, 1]$ open in $[0, 1]$ and
containing $t$ and a map $s: I \to B$ so that $\pi \circ s = id_I$.
We want to patch these local sections into a global section.

Since $[0, 1]$ is compact, we can cover it by finitely many intervals
$I_j$ so that on each $I_j$ there is a section $s_j: I_j \to B$.  Let
us now assume for simplicity that there are only two intervals: $I_0 =
[0, 2/3)$ and $I_1 = (1/3, 1]$.  The case of more than two intervals
will be left as an exercise to the reader.  Thus we have two sections
$s_0: [0, 2/3) \to B$ and $s_1: (1/3, 1] \to B$.  Since $G$ acts freely on
$\pi \inv ((0, 1))$ the map $\varphi : (1/3, 2/3) \times G \to \pi\inv
((1/3, 2/3))$ given by $\varphi (t, g) = g\cdot s_0 (t)$ is a
$G$-equivariant diffeomorphism (where $G$ acts on the product $(1/3,
2/3) \times G$ by multiplication on the second factor).  We have
$\varphi\inv \circ s_0 (t) = (t, 1)$ and $\varphi\inv \circ s_t (t) =
(t, g(t))$ for some curve $g: (1/3, 2/3) \to G$.  Since $(1/3, 2/3)$ is
simply connected we may lift $g$ to a curve $\gamma (t) $ the
universal cover $\exp : \fg\to G$ of $G$.  Choose a smooth function
$\rho: (1/3, 2/3) \to [0, 1]$ with $\rho (t) \equiv 0$ for $t$ near
$1/3$ and $\rho (t) \equiv 1$ for $t$ near $2/3$.  Now consider the
curve $a: (1/3, 2/3) \to G$ given by $a (t) = \exp ( \rho (t) \gamma
(t))$.  The map $s _2(t) = \varphi (t, a (t))$ is a local section of
$\pi :B \to [0,1]$ which agrees with $s_0$ near $1/3$ and with $s_1$
near $2/3$.  Thus we can define a global section $s: [0, 1] \to B$ by
$$
s(t) =
\begin{cases}
 s_0 (t ) & \, t\in [0, 1/3],\\
s_2 (t) & \, t\in [1/3, 2/3],\\
s_3 (t) & \, t\in [1/3, 0].
\end{cases}
$$ 
Now that we have a global section of $\pi : B\to [0, 1]$ we can
define a continuous map $\tilde{f} : [0, 1] \times G \to B$, $\tilde
{f}(t,g) = g\cdot s(t)$.  The map is onto; it descends to a bijective
continuous map $f: ([0, 1]\times G)/\!\sim \,\to B$ where $\sim$ is the
equivalence relation in the statement of the Lemma.  Since $([0,
1]\times G)/\sim$ is compact, the map $f$ is a homeomorphism.
\end{proof}

\begin{exercise}
Show that if the groups $K_1$ and $K_2$ in the statement of the Lemma
are the same, then $B$ is $S^1 \times S^2$.
\end{exercise}

\begin{exercise} (This exercise is considerably harder than the one above.)
Show that if the groups $K_1$ and $K_2$ are different, then $B$ is the
quotient of $S^3$ by a finite cyclic group.
\end{exercise}

\begin{lemma}\label{lemma-coh-lens}  Let $G= \bbT^2$ and let 
$K_1, K_2 \subset G$ be two closed subgroups isomorphic to $S^1$.  Let
$B$ be the topological space $([0, 1] \times G)/\sim$ where $(0, g)
\sim (0, ag)$ for all $g\in G$, $a \in K_1$ and $(1, g) \sim (1, ag) $
for all $g\in G$ and $a\in K_2$.  In other words $B$ is obtained from
the manifold with boundary $[0,1] \times G$ by collapsing circles in
the two components of the boundary by the respective actions of two
 circle subgroups.

Then either $H^1 (B, \Z) = \Z = H^2(B, \Z)$ or $H^1 (B, \Z) = 0$ and
$H^2 (B, \Z)$ is a finite group.  In particular, $B$ cannot be
homeomorphic to the 3-torus $\bbT^3$.
\end{lemma}

\begin{proof}
Recall that $H^1 (G, \Z)$ is isomorphic to the weight lattice $\Z_G^*$
and that the isomorphism is given as follows:  A weight $\nu \in
\Z_G^*$ defines a character $\chi_\nu : G \to S^1$ by $\chi_\nu ( \exp
(X)) = e^{2\pi i \nu (X)}$; the class $\chi_\nu ^* [d\theta]$ is the
element in $H^1 (G, \Z)$ corresponding to $\nu$.  Here $d\theta$ is
the obvious 1-form on $S^1$.

Consequently if $G= \bbT^2$ and $K_j \subset G$ is a circle subgroup,
then $\pi_j : G \to G/K_j \simeq S^1$ is a character and hence the
weight $\nu_j = (d \pi_j)_1 $ defines an element of $H^1 (G, \Z)$.
Thus if we identify $H^1 (G/K_j, \Z)$ with $\Z$ and $H^1 (G, \Z)$ with
$\Z_G^*$, then the map $H^1 (G/K_j, \Z) \to H^1 (G, \Z)$ becomes the
map $\Z\ni n\mapsto n\nu_j \in \Z_G^*$.

The sets $U = ([0, 2/3) \times G)/\sim $ and $V = ((1/3, 1] \times
G)/\sim$ are two open subsets of $B$.  We have $B = U \cup V$, $U\cap
V = (1/3, 2/3)\times G$ is homotopy equivalent to $G$, $U$ is homotopy
equivalent to $G/K_1$, $V$ is homotopy equivalent to $G/K_2$ and the
inclusion maps $U\cap V \hookrightarrow U$, $U\cap V \hookrightarrow
V$ are homotopy equivalent to projections $\pi_1: G \to G/K_1$, $\pi_2
:G \to G/K_2$ respectively.  Hence under the above identifications of
$H^1 (U)$ and $H^1(V)$ with $\Z$, the inclusions $U\cap V \to U$, $U\cap V
\to V$ induce the maps $\Z\ni n\mapsto n\nu_j \in \Z_G^*$, $j=1, 2$,
respectively.

We now apply the Mayer-Vietoris sequence to compute the integral
cohomology of $B$.  We have: $0 \to H^0(B) \to H^0(U) \oplus H^0
(V) \to H^0 (G)\stackrel{\delta}{ \to} H^1(B) \to H^1(U) \oplus H^1(V)
\to H^1 (G)
\stackrel{\delta}{ \to} H^2 (B)\to H^2(U) \oplus H^2(V) \to H^2 (G)
\stackrel{\delta}{ \to}H^3 (B) \to 0$.  Clearly the map $H^0(U) \oplus
H^0 (V) \to H^0 (G) $ is onto.  Given the identifications above the
map $\varphi :H^1(U) \oplus H^1(V) \to H^1 (G)$ becomes $\Z \oplus \Z\ni
(n, m) \mapsto n\nu_1 + m \nu_2 \in \Z_G^*$.   We therefore have
$0 \to H^1 (B) \to \Z \oplus \Z \stackrel{\varphi}{\to } \Z_G^*
\stackrel{\delta}{ \to} H^2 (B) \to 0\oplus 0 \to H^2 (G) 
\stackrel{\delta}{ \to} H^3 (B) \to 0$.
We conclude that 
\begin{itemize}
\item $H^2 (B ) = \Z_G^* / (\Z \nu_1 + \Z \nu _2)$

\item $H^1 (B) = \{(n, m) \in \Z^2 \mid n \nu_1 + m \nu_2 = 0\}$.

\end{itemize}
Since $\nu_1, \nu_2$ are differentials of projections onto quotients
by circle subgroups, either $\nu_1$ and $\nu_2$ are independent over
$\Z$ or $\nu_1 = \pm \nu_2$.  In the first case $H^1 (B) = 0$ and $H^2
(B) $ is a finite abelian group.  In the second case $H^1 (B) = \Z$
and $H^2 (B) = \Z_G^* /\Z\nu_1 = \Z$.
\end{proof}
This finishes the proof of Theorem~\ref{main-thm'} in the case that
$\dim B = 3$.

\chapter{\protect Proof of Theorem~\ref{main-thm'}, part II: 
uniqueness of symplectic toric manifolds}

In this section we sketch a proof of

\begin{theorem}\label{main-thm2}
Let $(B, \xi = \ker \alpha, \Psi: \xi^\circ_+ \to \fg^*)$ be a compact
connected contact toric $G$-manifold.

Suppose the dimension $k$ of the maximal linear subspace of the moment
cone $C(\Psi) = \Psi (\xi^\circ_+) \cup \{0\}$ satisfies $0 < k < \dim
G$.  Then $B$ is homotopy equivalent to the product of a $k$-torus
with a sphere.  In particular $B$ is not the co-sphere bundle of the
$n$-torus, $n = \frac{1}{2} (\dim B + 1) = \dim G$.
\end{theorem}

The main idea of the proof is simple.  We will first argue that there
is an action of $G$ on the symplectic manifold 
$$
M = T^* \bbT^k \times
\C^l \smallsetminus 0 = \{ (q,p,z) \in \bbT^k \times (\R^k)^* \times
\C^l = T^* \bbT^k \times \C^l \mid (p, z) \not = (0,0)\},
$$ 
$l = \frac{1}{2} (\dim B + 1 - k)>0$, with moment map $\tilde {\Phi} : M
\to \fg^*$ such that
\begin{equation} \label{eq***}
 \tilde{\Phi }(M) = \Psi (\xi^\circ _+).
\end{equation}
We then argue that (\ref{eq***}) implies that $M$ is $G$-equivariantly
symplectomorphic to $\xi^\circ_+$.  Note that $M$ is homotopy
equivalent to $\bbT^k \times S^{k+ 2l -1}$.

We start with the definition of a  symplectic slice representation
(c.f. proof of Lemma~\ref{lemma-noG}), which is essential for
understanding the local structure of symplectic toric manifolds.

\begin{definition}
Let $(M, \omega)$ be a symplectic manifold with a Hamiltonian action
of a torus $G$.  Then an orbit $G\cdot m$ is an isotropic submanifold
of $(M, \omega)$.\footnote{For a proof of this easy fact see, for
example, \cite{GSbook}. }  The {\bf symplectic slice representation}
at $m$ is the representation of the isotropy group $G_m$ of $m$ on the
symplectic vector space $V:= T_m (G\cdot m)^\omega /T_m (G\cdot m)$.
Here, as usual, $ T_m (G\cdot m)^\omega $ denotes the symplectic
perpendicular to $T_m (G\cdot m)$ in the symplectic vector space $(T_m
M, \omega_m)$.
\end{definition}

The equivariant isotropic embedding theorem (see for example
\cite{GSbook}, Theorem~39.1) asserts that a neighborhood of an orbit
$G\cdot m$ is determined (up to equivariant symplectomorphisms) by the
symplectic slice representation at $m$. In fact, the topological
normal bundle of the orbit $G\cdot m$ in $M$ is $G\times _H (\fg/\fh
\times V)$ where $H$ is the isotropy group of $m$ and $\fh$ is its Lie
algebra ({\em op.\ cit.}).

\begin{remark}
In the case of symplectic toric manifolds the dimension of the
symplectic slice $V$ at $m$ is twice the dimension of the isotropy
group $G_m$.  Hence by Lemma~\ref{Del-lemma} the group $G_m$ is
connected and the representation of $G_m$ on $V$ is determined by a
set of weights $\{ \nu_i\}$ which forms a basis of the weight lattice
of $G_m$.  Note that the image of $V$ under the moment map $\Phi_V : V
\to \fg_m^*$ defined by the representation is the cone
$$
\Phi_V (V)  = \{ \sum a_i \nu_i \mid a_i \geq 0\}.
$$ 
In particular the edges of the cone are spanned by the weights.
Alternatively, the isotropy group $G_m$ is isomorphic to $\bbT^l$ for
some $l$ and the slice representation $\rho: G_m \to \text{Sp}(V)$ is
isomorphic to the standard representation of $\bbT^l$ on
$\C^l$.\footnote{ The standard representation of $\bbT^l = \{
(\lambda_1, \ldots, \lambda_l) \in \C^l \mid |\lambda _j | = 1\}$ on
$\C^l$ is given by $(\lambda_1, \ldots, \lambda_l) \cdot (z_1, \ldots,
z_l) = (\lambda_1 z_1, \ldots, \lambda_l z_l)$.}
\end{remark}

With a little more work one can prove the following two propositions
(their proofs can be found, for example, in \cite{D}).
\begin{proposition} \labell{prop1}
Let $(M, \omega, \Phi : M\to \fg^*)$ be a symplectic toric manifold,
$m\in M$ a point and $U$ a neighborhood of $G\cdot m$ in $M$. Then for
a sufficiently small ball $\cO$ about $\eta = \Phi (m)$ in $\fg^*$ the
set $\Phi (U) \cap \cO$ determines the symplectic slice representation
at $m$ and hence a small $G$-invariant neighborhood of $G\cdot m$ in
$M$.
\end{proposition}
\begin{proposition}\label{prop2}
Let $(M, \omega, \Phi : M\to \fg^*)$ be a symplectic toric manifold,
$m\in M$ a point, $\eta = \Phi (m)$ and $G_m$ the isotropy group of
$m$.  Identify $G_m$ with the standard torus $\bbT^d$, $d = \dim G_m$,
and extend it to an identification of $G$ with $\bbT^d \times \bbT^c$,
\footnote{Here we use Remark~\ref{splitting-tori}.} 
$c = \dim G - \dim G_m$.

A $G$-invariant neighborhood $U$ of $G\cdot m$ in $M$ is equivariantly
symplectomorphic to a neighborhood of $\bbT^c \times \{(0, 0)\}$ in
$T^* \bbT^c\times \C^d \simeq \bbT^c \times (\R^c)^* \times \C^d$.  Hence
$$
\left\{ t (\mu - \eta) + \eta \mid t\geq 0, \, \mu \in \Phi (U)\right\} =
\tilde{\Phi} (T^* \bbT^c \times \C^d),
$$ where $\tilde{\Phi} : T^* \bbT^c \times \C^d \to (\R^c)^* \times
(\R^d)^* \simeq \fg^*$ is the moment map for the ``obvious'' action of
$G = \bbT^c \times \bbT^d$ on $T^* \bbT^c \times \C^d$.\footnote{The
``obvious'' action of $\bbT^c$ on $T^* \bbT^c$ is the lift of left
multiplication.}
\end{proposition}

\begin{corollary}
Let $(B, \xi = \ker \alpha, \Psi: \xi_+^\circ)$ be a compact connected
contact toric $G$-manifold.  Suppose the dimension $k$ of the
maximal linear subspace of the moment cone $C(\Psi) =\Psi (\xi^\circ_+) \cup
\{0\}$ satisfies $0 < k < \dim G$.

Then there is an identification of $G$ with $\bbT^k \times \bbT^l$, 
$l = \dim G - k$, so that
$$
\Psi (\xi^\circ_+) \cup \{0\} = \tilde{\Phi} (T^* \bbT^k \times \C^l)
$$ 
where $\tilde{\Phi}:T^* \bbT^k \times \C^l \to (\R^k)^* \times
(\R^l)^* \simeq \fg^*$ is the moment map for the obvious action of
$\bbT^k \times \bbT^l $ on $T^* \bbT^k \times \C^l$.
\end{corollary}

\begin{proof}
If $C$ is a cone in $\fg^*$ whose maximal linear subspace is $P$ and
if $U$ is a neighborhood in $\fg^*$ of a point $\eta\in P$ then $C$
equals the cone on $U\cap C$ with the vertex at $\eta$:
$$
C = \{ t (\mu - \eta) + \eta \mid t\geq 0, \mu \in U\cap C\}.
$$

 Since the $B$ is contact toric $\Psi (\xi^\circ_+) $ does
not contain the origin (c.f. Lemma~\ref{lemma-noG}).  Since $B$ is
compact and $\Psi: \xi^\circ_+ \to \fg^*$ is homogeneous, $\Psi :
\xi^\circ_+ \to \fg^* \smallsetminus \{0\}$ is proper. Hence for any 
$\eta \in \Psi (\xi^\circ_+)$ and any neighborhood $U$ of $\Psi\inv
(\eta)$, $\Psi (U)$ is a neighborhood of $\eta$ in $\Psi
(\xi_+^\circ)$.  

Moreover, since the fibers of $\Psi$ are $G$-orbits (by
Theorem~\ref{thm-BM}), for any neighborhood $U$ of an orbit $G\cdot m$
the set $\Psi (U)$ is a neighborhood of $\eta = \Psi (m)$ in $\Psi
(\xi^\circ_+)$.  In particular it contains a set of the form $\cO \cap
\Psi (\xi^\circ_+)$ where $\cO \subset \fg^*$ is a small ball about
$\eta$.

Let $P\subset C(\Psi)$ be the maximal linear subspace and let $0\not =
\eta \in P$.  Then, as noted at the beginning of the proof, for any
ball $\cO\subset \fg^*$ about $\eta$ we have
\begin{equation*}
\begin{split}
C(\Psi) & = 
\left\{ t (\mu - \eta) + \eta \mid t\geq 0, \, 
\mu \in \cO \cap \Psi (\xi^\circ_+)\right\}\\
&= \left\{ t (\mu - \eta) + \eta \mid t\geq 0, \, \mu \in \Psi (U)\right\}.\\
\end{split}
\end{equation*}
On the other hand, it follows from Proposition~\ref{prop2} that
$$
\left\{ t (\mu - \eta) + \eta \mid t\geq 0, \, \mu \in \Psi (U)\right\} =
\tilde{\Phi} (T^* \bbT^c \times \C^d) \simeq \R^c \times (\R_{\geq 0})^d
$$ for some integers $c$ and $d$.  It follows that $c=k$, $d = \dim G
- k = l$ and that
$$
C(\Psi) = \tilde{\Phi} (T^* \bbT^k \times \C^l).
$$ 
\end{proof}

Note that since $\tilde{\Phi}\inv (0) = \bbT^k \times \{(0,0)\}$, we get
$$
\Psi (\xi^\circ_+) = \tilde{\Phi}(M),
$$
where $M = (T^* \bbT^k \times \C^l )\smallsetminus  (\bbT^k \times \{(0,0)\})$.

\begin{definition}
A symplectic toric $G$-manifold $(M, \omega, \Phi: M \to \fg^*)$ is
{\bf good} (for the purposes of these lectures) if
\begin{enumerate}
\item the fibers of $\Phi : M \to \fg^*$ are connected,
\item the set $\Phi (M)$ is contractible and 
\item there is an open set $U\subset \fg^*$ such that $\Phi (M) \subset U$ 
and the map $\Phi :M \to U$ is proper.
\end{enumerate}
\end{definition} 
The definition is designed in such a way that it includes compact
symplectic toric manifolds,  symplectizations of compact contact
toric manifolds of dimension bigger than 3 and the manifolds of the form 
$M = (T^* \bbT^k \times\C^l )\smallsetminus  (\bbT^k \times \{(0,0)\})$, 
$0<k, l$.  As a consequence of the
definition and of Proposition~\ref{prop1} we have

\begin{lemma}\label{lemma*1} 
Let $(M, \omega, \Phi : M\to \fg^*)$ be a good symplectic toric
manifold.  For any point $\eta \in \Phi (M)$ and any sufficiently
small ball $\cO$ centered at $\eta$, the set $\cO \cap \Phi (M)$
determines the symplectic toric manifold $(\Phi\inv (\cO),
\omega|_{\Phi\inv (\cO)},
\Phi|_{\Phi\inv (\cO)})$ (up to a $G$-equivariant symplectomorphism).
\end{lemma}

\begin{definition}
Two symplectic toric $G$-manifolds $(M, \omega, \Phi: M\to \fg^*)$ and
$(M', \omega', \Phi': M'\to \fg^*)$ are {\bf isomorphic} if there is a
$G$-equivariant diffeomorphism $\sigma : M\to M'$ such that $\sigma^*
\omega' = \omega$ and $\sigma^* \Phi' = \Phi$.

We denote by $\text{Iso} (M, \omega, \Phi) = \text{Iso} (M)$ the group
of isomorphisms of a symplectic toric manifold $(M, \omega, \Phi)$.
\end{definition}
Note that the last condition on $\sigma$ is almost redundant --- if
$\sigma$ is symplectic and $G$-equivariant then $\sigma^* \Phi' = \Phi
+ c$ for some constant vector $c\in \fg^*$.  We impose the last
condition for technical convenience.
\begin{definition}
Two symplectic toric $G$-manifolds $(M, \omega, \Phi: M\to \fg^*)$ and
$(M', \omega', \Phi': M'\to \fg^*)$ are {\bf locally isomorphic} over
a set $\Delta \subset \fg^*$ if
\begin{enumerate}
\item $\Phi (M) = \Delta = \Phi' (M')$ and
\item for any $\eta \in \Delta$ and any sufficiently small ball 
$\cO\subset \fg^*$ centered at $\eta$ the symplectic toric manifolds 
$(\Phi\inv (\cO), \omega|_{\Phi\inv (\cO)},
\Phi|_{\Phi\inv (\cO)})$ and $((\Phi')\inv (\cO), \omega'|_{(\Phi')\inv (\cO)},
\Phi'|_{(\Phi')\inv (\cO)})$ are isomorphic.
\end{enumerate}
\end{definition}
Given the above definitions we see that Lemma~\ref{lemma*1} implies

\begin{lemma}
Suppose $(M, \omega, \Phi)$ and $(M', \omega', \Phi')$ are two good
symplectic toric $G$-manifolds with $\Phi (M) = \Phi'(M')$.  Then $(M,
\omega, \Phi)$ and $(M', \omega', \Phi')$ are locally isomorphic over
$ \Delta = \Phi (M)$.
\end{lemma}

Therefore the proof of Theorem~\ref{main-thm2} reduces to 
\begin{proposition} \labell{prop101}
Any good symplectic toric $G$-manifold locally isomorphic to a given
 symplectic toric $G$-manifold $(M, \omega, \Phi)$ is actually
isomorphic to $(M, \omega, \Phi)$.
\end{proposition}
The proposition could be stated for a larger class of symplectic toric
manifolds.  We leave it to the reader to find the most general form of
the statement (and prove it).  The rest of the section is occupied
with a proof of the proposition.

Suppose a symplectic toric $G$-manifold $(M', \omega', \Phi')$ is
locally isomorphic to $(M, \omega, \Phi)$.  Then for any $\eta \in
\Phi (M)$ and for any sufficiently small ball $\cO\subset \fg^*$ about
$\eta$ the symplectic toric manifold $\Phi\inv (\cO)$ is isomorphic to
$(\Phi')\inv (\cO)$.  Choose a locally finite cover $\{\cO_i\}_{i\in I}$
of $\Phi(M)$ by such balls.  Then for each index $i$ we have an
isomorphism
$$
f_i:  \Phi\inv (\cO_i)\to (\Phi')\inv (\cO_i) .
$$
Let $\cO_{ij} = \cO_i \cap \cO_j$.  Define 
$$
g_{ij} : \Phi\inv (\cO_{ij}) \to  \Phi\inv (\cO_{ij}), \quad 
g_{ij} = f_i\inv \circ f_j
$$
(to keep the notation manageable we wrote $f_i$ for the restriction
$f_i|_{\cO_{ij}}$ etc.; we will continue to omit restrictions in the
rest of the section).  It is easy to see that
\begin{equation}\labell{eq-star}
g_{ii} = id, \, g_{ij}\circ g_{ji} = id, \,\text{and }\, 
g_{ij} \circ g_{jk} \circ g_{ki} = id
\end{equation}
wherever these equations make sense.

The data $\{\cO_i\}$, $\{g_{ij}\}$ and $(M, \omega, \Phi)$ allow us to
reconstruct $(M', \omega', \Phi')$.  Indeed, let $\tilde{M} =
\bigsqcup \Phi\inv (\cO_i)$.  Define a relation $\sim$ on $\tilde{M}$ by 
$\Phi\inv (\cO_i)\ni x_i \sim x_j \in \Phi\inv (\cO_j)$ iff $x_j =
g_{ij} (x_i)$.  Equations (\ref{eq-star}) imply that $\sim$ is an
equivalence relation. It follows that $\tilde{M}/\sim$ is a symplectic
toric manifold. Moreover the map $\tilde{F} : \tilde{M} \to M'$ defined
by $\tilde{F}|_{\Phi\inv (\cO_i)} = f_i$ descends to a well
defined map $F: \tilde{M}/\sim \to M'$.  The map $F$ is an
isomorphism.

Suppose $\tilde{f}_i : \Phi\inv (\cO_i) \to (\Phi') \inv (\cO_i)$ is
another collection of isomorphisms.  Let $\tilde{g}_{ij} =
(\tilde{f}_i)\inv \circ \tilde{f}_j$.  Clearly
\begin{equation}\labell{eq6.3}
g_{ij} =  h_i \circ \tilde{g}_{ij} \circ h_j\inv
\end{equation}
where $h_i = f_i \inv \circ \tilde{f}_i \in \text{Iso} (\Phi\inv (\cO_i))$.

We also get a different set of data if we choose a different
cover. However, all these sets of data define one object --- a class
in the first \v{C}ech cohomology with coefficients in a certain
sheaf.  
Let us now review the notions of sheaves and \v{C}ech cohomology.
There are many good references for this material.  I will be following
\cite{WW}.

\begin{definition}
A {\bf sheaf} of groups $\cS$ on a topological space $X$ is an assignment
$$
\cS:\{ \text{opens sets in }X\} \to \text{groups}, \,\, U\mapsto \cS(U)
$$
satisfying two conditions:
\begin{enumerate}
\item For a pair of open sets $U\subset W$ in $X$ there is a restriction map
$\rho^W_U: \cS (W) \to \cS(U)$ such that for any three sets $U\subset
W \subset V$ of $X$ 
$$
\rho^W_U \circ \rho^V_W = \rho^V_U.
$$
Elements of $\cS(U)$ are called {\bf sections}.  Given a section
$\varphi \in \cS (W)$ we write $\varphi |_U$ for $\rho^W_U (\varphi)$.

\item Given an open cover $\{U_i\}$ of an open set $U$ 
(so that $U = \bigcup U_i$) and a collection of sections 
$\varphi_i \in \cS(U_i)$ such that
$$
\varphi_i |_{U_i \cap U_j} = \varphi_j |_{U_i \cap U_j}
$$
for all indices $i,j$ there is a unique section $\varphi \in \cS (U)$
such that
$$
\varphi |_{U_i} = \varphi_i
$$
for all $i$.
\end{enumerate}
The sheaf $\cS$ is {\bf abelian} if $\cS(U)$ is an
abelian group for all open sets $U$.
\end{definition} 
Three examples of sheaves will be important to us.  Check that they
are indeed sheaves.
\begin{example}\label{ex-iso}
Let $(M, \omega, \Phi)$ be a good symplectic toric $G$-manifold.  The
assignment $\text{Iso}: U\mapsto \text{Iso}(U)$ ($U\subset \Phi(M)$ open) is
a sheaf on $\Phi (M)$.  The group operation is composition.
\end{example}
\begin{example}[Locally constant sheaf]
Let $H$ be a group and $X$ a topological space.  The assignment that
associates the group $H$ to every open connected subset of $X$ is a
sheaf, called a locally constant sheaf.  It is denoted by
$\underline{H}$. Thus $\underline{H} (U) = H$ for every connected open
set $U \subset X$.  The group operation is the multiplication in $H$.
\end{example}
\begin{example}\label{ex-c}
Let $(M, \omega, \Phi)$ be a good symplectic toric $G$-manifold.  
Define a sheaf $\cC$ on $\Phi(M)$ by 
$$
\cC (U) = C^\infty (\Phi\inv (U))^G, \quad 
\text{$G$-invariant smooth functions on $\Phi\inv (U)$}.
$$
The group operation is addition of functions.
\end{example}

\begin{definition}
Let $\cS_1, \cS_2$ be sheaves on a topological space $X$.  A {\bf map
of sheaves} $\tau: \cS _1 \to \cS_2$ is a family of group
homomorphisms $$
\tau_U : \cS_1 (U) \to \cS_2 (U), \quad U\subset X \text{ open} 
$$
compatible with the restrictions:
$$
(\rho_2)^W_U \circ f_W = f_U \circ (\rho_1)^W_U 
\quad \text{ for all pairs } U\subset W \text{ of open sets in } X.
$$
\end{definition}

\begin{example}
Consider the sheaves $\text{Iso}$ and $\cC$ defined in
Examples~\ref{ex-iso} and \ref{ex-c} above.  A section $f\in \cC(U)$
is a $G$-invariant function on $\Phi\inv (U)$.  Its time $t$ flow
$\varphi_t^f$ preserves the fibers of $\Phi$ and is a $G$-equivariant
symplectomorphism of $\Phi\inv (U)$.  Hence $\varphi^f_1$ is a section
of $\text{Iso}(U)$.  This gives us for each open set $U\subset
\Phi(M)$ a map $\tau_U : \cC (U) \to \text{Iso} (U)$, $\tau_U (f) =
\varphi^f_1$.   Moreover,  any two functions $f_1, f_2 \in \cC(U)$ Poisson 
commute [prove this].  Hence $\varphi^{f_1}_t \circ \varphi^{f_2}_t =
\varphi^{f_1 + f_2}_t$ and therefore $\tau_U$ is a group homomorphism.
Thus we get a map of sheaves $\tau :\cC \to \text{Iso}$.
\end{example}

Given a map of sheaves $\tau : \cS_1 \to \cS_2$ one can define the
sheaves the kernel and image sheaves $\ker \tau $ and $\text{im}
\tau$: for an open set $U$ $(\ker \tau) (U) = \ker \tau_U$, 
$(\im \tau)(U) = \im \tau U$.  Hence it makes sense to say that a map
of sheaves is onto and more generally talk about exact sequences of
sheaves.\footnote{ Warning: the map $\tau : S_1 \to S_2$ being onto
{\em does not} mean that $\tau_U$ is onto for every open set $U$. See
\cite{WW} or any other good book on sheaves for more details.}

\begin{proposition}\label{prop-}
Let $(M, \omega, \Phi)$ be a good symplectic toric manifold.  The
map of sheaves $\tau : \cC \to \text{\rm Iso}$ defined above is  onto.  
Hence $\text{\rm Iso}$ is an abelian sheaf.

The kernel of $\tau $ is the locally constant sheaf $\underline{\R
\times \Z_G }$. We thus have a short exact sequence of abelian
sheaves: 
$$ 
0\to \underline{\R \times \Z_G} \to \cC \to \text{ \rm
Iso} \to 0.  
$$
\end{proposition}
The proposition is due to Boucetta and Molino \cite{BoM}.  See also
\cite{LT}.

\section{ \v{C}ech cohomology}.

\noindent
In this subsection we ``review'' the notion of \v Cech cohomology with
coefficients in an abelian sheaf.  There are many good references,
such as \cite{WW}, for the nontrivial facts that we list below
without proofs.

Let $X$ be a topological space $\{U_i\}$ an open locally finite cover
of $X$ and $\cS$ an abelian sheaf on $X$. A 0 \v Cech cochain is a
function that assigns to each index $i$ an element $f_i$ of $\cS
(U_i)$, i.e., the group of $0$-cochains $C^0 (\{U_i\}, \cS)$ is the
product $\prod \cS(U_i)$. A 1
\v Cech cochain assigns to an ordered pair of indices $ij$ an element
$g_{ij}$ of $\cS(U_{ij})$ where $U_{ij} = U_i \cap U_j$.  Moreover we
require that $g_{ij} = - g_{ji}$ (we now think of the groups $\cS(U)$
additively).  More generally a $p$-cochain assigns to an ordered $p+1$
tuple of indices $i_0\ldots i_p$ an element $s_{i_0\ldots i_p}\in \cS
(U_{i_0\ldots i_p})$ where $U_{i_0\ldots i_p} = U_{i_0} \cap \ldots
\cap U_{i_p}$ and $s_{i_0\ldots i_p}$ is skew-symmetric in the
indices.  The coboundary operator $\delta: C^p (\{U_i\}, \cS) \to
C^{p+1}(\{U_i\}, \cS)$ is defined by 
$$ 
(\delta s)_{i_0\ldots i_{p+1}}
= \sum (-1)^j s_{i_0\ldots \hat{i}_j \ldots i_{p+1}} 
$$ 
where $ \hat{i}_j$
means that the index is omitted, and where we omitted writing the
restrictions of the terms on the right hand side to $U_{i_0\ldots
i_{p+1}}$. One proves that $\delta^2 = 0$. The cohomology of the
complex $( C^p (\{U_i\}, \cS), \delta)$ denoted by $\check H ^*
(\{U_i\}, \cS)$ is called the \v Cech cohomology of the cover
$\{U_i\}$ with coefficients in the sheaf $\cS$.

Given a refinement $\{V_j\}$ of the cover $\{U_i\}$ the restrictions
give rise to a chain map $C^p (\{U_i\}, \cS) \to C^p (\{V_j\}, \cS)$,
which in turn gives rise to a map in cohomology $\check H ^* (\{U_i\},
\cS) \to \check H^* (\{V_j\}, \cS)$.  Taking the direct limit over all 
locally finite covers we get a well-defined cohomology group
$$
\check H^* (X, \cS) =\lim _\to \check H ^* (\{U_i\}, \cS),
$$ 
the \v Cech cohomology of $X$ with coefficients in the sheaf $\cS$.

Now let $(M, \omega, \Phi: M\to \fg^*)$ be a good symplectic toric
$G$-manifold and $\{\cO_i\}$ a locally finite cover of $\Phi(M)$ by
sufficiently small balls.  If $\{g_{ij}\} \in C^1 (\{\cO_i\},
\text{Iso})$ is a 1-cochain, then $\delta (\{g_{**}\}) = 0$ means that
for all triples of indices $ijk$ we have 
$$ 
0 = \delta
(\{g_{**}\})_{ijk} = -g_{jk} + g_{ik} + g_{ij}, 
$$ 
which is (\ref{eq-star}) in additive notation (where on the right hand
side we omitted the restrictions to $U_{ijk}$).  Similarly if
$\{g_{ij}\},\{\tilde{g}_{ij}\} \in C^1 (\{\cO_i\},
\text{Iso})$ are two  1-cochains that differ by $\delta(\{h_*\})$ 
for some 0-cochain $\{h_i\}$ then 
$$
\tilde{g}_{ij} - g_{ij} = - h_i + h_j
$$ 
hence 
$$ 
g_{ij} = h_i + \tilde{g}_{ij} - h_j, 
$$ 
which is (\ref{eq6.3}) in additive notation.  Thus the discussion
above shows that to every element of $\check H^1 (\{\cO_i\},
\text{Iso})$ there corresponds a good symplectic toric manifold $(M',
\omega', \Phi' :M\to \fg^*)$ locally isomorphic to 
$(M, \omega, \Phi: M\to \fg^*)$.  More generally one can check that
there is a one-to-one correspondence between cohomology classes in
$\check H^1 (\Phi (M), \text{Iso})$ and isomorphism classes of good
symplectic toric manifolds locally isomorphic to $(M, \omega, \Phi:
M\to \fg^*)$.  Thus to complete the proof of Proposition~\ref{prop101}
(and thereby Theorem~\ref{main-thm2}) it remains to show that the
group $\check H^1 (\Phi (M), \text{Iso})$ is trivial for any good
symplectic toric manifold $M$.  For this we use
Proposition~\ref{prop-}, two properties of \v Cech cohomology and a
property of the sheaf $\cC$ defined in Example~\ref{ex-c}.
The first property of \v Cech cohomology that we need is 
\begin{theorem}\label{thm*}
A short exact sequence of abelian sheaves $0 \to \cS_1 \to \cS_2 \to
\cS_3 \to 0$ on a space $X$ induces a long exact sequence in \v Cech
cohomology 
$$
\cdots \to \check H^p (X, \cS_1) \to\check H^p (X, \cS_2) \to
\check H^p (X, \cS_3) \stackrel{\delta}{\to}\check H^{p+1} (X, \cS1) \to \cdots
$$
\end{theorem}
\noindent
The second property that we will use is 
\begin{theorem}\label{thm-ch-sing}
Let $X$ be a simply connected topological space and $H$ an abelian
group.  The \v Cech cohomology $\check H ^* (X, \underline{H})$ of $X$
with coefficients in the locally  constant sheaf $\underline{H}$ is isomorphic
to the singular cohomology $H^* (X, H)$ of $X$ with coefficients in
the abelian group $H$.
\end{theorem}
\noindent
We will use the following property of the sheaf $\cC$ (cf.\ \cite{LT},
Proposition~7.3)
\begin{lemma}\label{lem-3*}
The sheaf $\cC$ defined in Example~\ref{ex-c} is acyclic, that is,\\
$\check H^q (\Phi (M), \cC) = 0$ for all $q>0$.
\end{lemma}
Now putting Theorem~\ref{thm*}, Lemma~\ref{lem-3*} and
Proposition~\ref{prop-} together we see that if $(M, \omega, \Phi)$ is
a good symplectic toric $G$-manifold and Iso the sheaf defined in
Example~\ref{ex-iso} then the cohomology group $\check H^1 (\Phi (M),
\text{Iso})$ is isomorphic to $\check H^2 (\Phi (M), \underline{\R
\times \Z_G})$.  The latter group is isomorphic to the singular
cohomology group $ H^2 (\Phi (M), \Z_G \times \R)$ by
Theorem~\ref{thm-ch-sing}.  But $\Phi(M)$ is contractible, so $ H^2
(\Phi (M), \Z_G \times \R) = 0$.  Therefore $\check H^1 (\Phi (M),
\text{Iso}) =0$, which proves Proposition~\ref{prop101} and thereby
Theorem~\ref{main-thm2}.

\chapter{\protect Proof of Theorem~\ref{main-thm'}, part III: Morse theory 
on orbifolds} The goal of this section is to prove
\begin{theorem}\labell{thm-}
Let $(B, \xi = \ker \alpha, \Psi: \xi^\circ_+ \to \fg^*)$ be a compact
connected contact toric $G$-manifold with $\dim B \geq 3$.  Suppose
there is a vector $X$ in the Lie algebra $\fg$ of $G$ such that the
function $\langle \Psi, X \rangle$ is strictly positive on $B$. Then
$\dim H^1 (B, \R) \leq 1$.  In particular $B$ is not the co-sphere
bundle $S^*G = \bbT^n \times S^{n-1}$, $n= \dim G = \frac{1}{2} (\dim
B +1) \geq 2$.
\end{theorem}
The proof of Theorem~\ref{thm-} above will complete our proof of
Theorem~\ref{main-thm'}.  Since Theorem~\ref{main-thm'} implies the
main result of the notes, Theorem~\ref{main}, this, in turn, will
finish the proof of the main result.
As was sketched out at the beginning of Chapter~4 our proof of
Theorem~\ref{thm-} has several steps.  The first one is a theorem
implicit in a paper of Boyer and Galicki \cite{BG}:
\begin{theorem}\labell{BGthm}
Let $(B, \xi = \ker \alpha, \Psi: \xi^\circ_+ \to \fg^*)$ be a compact
connected contact toric $G$-manifold with $\dim B \geq 3$.  Suppose
there is a vector $X$ in the Lie algebra $\fg$ of $G$ such that the
function $\langle \Psi, X \rangle$ is strictly positive on $B$. Then
there exists on $B$ a locally free circle action so that the quotient
$M = B/S^1$ is a (compact) symplectic toric orbifold.
\end{theorem}
The second step is the argument that if $M$ is a compact connected
symplectic toric orbifold then $H^q (M, \R) = 0$ for all {\em odd}
degrees $q$.  This step uses Morse theory on orbifolds.  

Let us now see why these two steps give us a proof of
Theorem~\ref{thm-}.  Consider the circle action produced by
Theorem~\ref{BGthm} and the corresponding $S^1$ orbit map $\pi : B \to
M$.  If the circle action is actually free, then $\pi$ is a circle
fibration and we have the Gysin sequence
\begin{equation}\label{eq-Gys}
0 = H^{-1} (M, \R) \to H^1 (M, \R) \to H^1 (B, \R) \to H^0 (M, \R) \to
H^2 (M, \R) \to \cdots .
\end{equation}
If the action of $S^1$ is locally free, the long exact sequence
(\ref{eq-Gys}) still exists.  The reason is that the Gysin sequence
arises from a collapse of the Leray-Serre spectral sequence for a
sphere bundle.  For locally free $S^1$ actions the orbit map $\pi:B\to
M$ is not a fibration, but the corresponding spectral sequences still
collapses if we use real coefficients.\footnote{ The reader not
familiar with spectral sequences may wish to take this claim on faith.
}  Now, since $H^1 (M, \R) = 0$ by the second step, it follows from
(\ref{eq-Gys}) that $\dim H^1 (B, \R)
\leq \dim H^0 (M, \R) = 1$, which proves Theorem~\ref{thm-}.

Here is an outline of the rest of the section.  We start with a review
of the notion of orbifold.  We recall how orbifolds arise as quotients
of manifolds by locally free actions of compact Lie groups.  We then
review the symplectic reduction theorem, in particular the fact that
generically symplectic reduced spaces are orbifolds and use it to
prove Theorem~\ref{BGthm}.  Next we define Morse functions on
orbifolds, and discuss the two fundamental results of Morse theory on
orbifolds. Finally we argue that a compact symplectic toric orbifold
$M$ has a Morse function with all indices even and use this to
conclude that the real cohomology of $M$ vanishes in odd degrees.

We now define orbifolds and related differential geometric notions.
For more details, see Satake \cite{Satake56, Satake57} and, for a more
modern point of view, Ruan \cite{R}.  The notion of orbifold was
introduced by Satake in 1956 under the name of V-manifold.  Orbifolds
are designed to generalize manifolds in the following sense: an
$n$-dimensional manifold is locally modeled on an open subset of $\R
^n$.  An $n$-dimensional orbifold is locally modeled on a quotient
$\tU /\Gamma$ where $\tU$ is an open subset of $\R^n$ and $\Gamma$ is
a finite group acting smoothly on $\tU$.  To give a precise definition
we need a few preliminary notions --- we need to define charts, atlases
and compatibility of atlases.

Let $U$ be a connected topological space, $\tU$ a connected
$n$-dimensional manifold and $\Gamma$ a finite group acting smoothly
on $\tU$.  An {\bf $n$-dimensional uniformizing chart}\footnote{Ruan
calls it a {\em uniformizing system.}} on $U$ is a triple $(\tU,
\Gamma, \varphi)$ where $\varphi: \tU \to U$ is a continuous map
inducing a homeomorphism between $\tU/\Gamma$ and $U$ (thus, in
particular, $\varphi$ is constant on the orbits of $\Gamma$).  We will
only consider charts where the set of points in $\tU$ fixed by
$\Gamma$ is either all of $\tU$ or is of codimension 2 or greater.
Note that we {\bf do not} require that $\Gamma$ acts effectively.

Two uniformizing charts $(\tU_1, \Gamma_1, \varphi_1)$ and $(\tU_2,
\Gamma_2, \varphi_2)$ of $U$ are {\bf isomorphic} if there is a
diffeomorphism $\psi: \tU_1 \to \tU_2$ and an isomorphism $\lambda :
\Gamma_1 \to \Gamma_2$ such that $\psi $ is $\lambda$-equivariant
(i.e., $\psi (g\cdot x) = \lambda (g) \cdot \psi (x)$ for all $g\in
\Gamma_1$, $x\in \tU_1$) and $\varphi_2 \circ \psi = \varphi_1$.
For example, fix $a\in \Gamma$.  Define $\psi :\tU \to \tU$ by $\psi
(x) = a \cdot x$.  Let $\lambda (g) = aga\inv$.  Then $(\psi,
\lambda): (\tU, \Gamma, \varphi) \to (\tU, \Gamma, \varphi)$ is an
isomorphism.

Let $\iota : U' \hookrightarrow U$ be a connected open subset of $U$.
We say that a uniformizing chart $(\tU', \Gamma', \varphi')$ is {\bf
induced from} a uniformizing chart $(\tU, \Gamma, \varphi)$ on $U$ if
there is a monomorphism $\lambda : \Gamma' \to \Gamma$ and a
$\lambda$-equivariant embedding $\psi: \tU' \to \tU$ such that $\iota
\circ \varphi' = \varphi \circ \psi$.  In this case $(\psi, \lambda) : 
(\tU', \Gamma', \varphi') \to (\tU, \Gamma, \varphi)$ is called an
{\bf injection}.

An {\bf orbifold atlas} on a Hausdorff topological space $M$ is an open cover $
\tcU$ of $M$ satisfying the following conditions:
\begin{enumerate}
\item  Each element $U$ of $\tcU$ is uniformized, say by 
$(\tU, \Gamma, \varphi)$.

\item If $U, U' \in \tcU$ and $U' \subset U$ then there is an injection 
$(\tU', \Gamma' , \varphi') \to (\tU, \Gamma, \varphi)$.

\item For any point $p\in U_1 \cap U_2$, $U_1, U_2 \in \tcU$, there is a 
connected open set $U_3 \in \tcU$ with $p\in U_3 \subset U_1 \cap U_2$ 
(and hence there are injections 
$(\tU_3, \Gamma_3 , \varphi_3) \to (\tU_1, \Gamma_1, \varphi_1)$ and
$(\tU_3, \Gamma_3 , \varphi_3) \to (\tU_2, \Gamma_2, \varphi_2)$).
\end{enumerate}

Suppose that $\tcV$ and $\tcU$ are two orbifold atlases on a space
$M$, that $\tcV$ is a refinement of $\tcU$ and that for every $V\in
\tcV$ and $U\in \tcU$ with $V\subset U$ we have an injection 
$(\tV, \Delta, \phi) \hookrightarrow (\tU, \Gamma, \varphi)$, where
$(\tV, \Delta, \phi)$ and $(\tU, \Gamma, \varphi)$ are the respective
uniformizing charts.  We then say that $\tcV$ and $\tcU$ are directly
equivalent orbifold atlases.  Now take the smallest equivalence
relation on the orbifold atlases on $M$ so that any two directly
equivalent atlases on $M$ are equivalent.  We now define an {\bf
orbifold} to be a Hausdorff topological space together with an
equivalence class of orbifold atlases.

Let $x$ be a point in an orbifold $M$, and let $(\tU,\Gamma, \varphi)$
be a uniformizing chart with $x\in\varphi (\tU)$.  The {\bf (orbifold)
structure group} of $x$ is the isotropy group of a point in the fiber
$\varphi\inv (x)$.  It is well-defined as an abstract group: if $x_1,
x_2 \in \varphi\inv (x)$ then the corresponding isotropy groups are
conjugate in $\Gamma$.

\begin{remark}\label{centered}
It is not hard to show that if $(\tU, \Gamma, \varphi)$ is a
uniformizing chart of $U$, then for any point $x\in U$ there is a
neighborhood $U'$ and a uniformizing chart $(\tU', \Gamma', \varphi')$
induced from $(\tU, \Gamma, \varphi)$ such that $(\varphi')\inv (x) $
is a single point $\tilde{x} $ (and hence $\Gamma'$ fixes $\tilde{x}$).
We will refer to $(\tU', \Gamma', \varphi')$ as a chart {\bf centered}
at $x$.
\end{remark}

Let $M$ be an orbifold with an atlas $\{\tU_i, \Gamma_i, \varphi_i\}$.
A {\bf smooth function} $f$ on $M$ is a collection of smooth
$\Gamma_i$-invariant functions $\tf_i$ on $\tU_i$ such that for any
injection $(\psi_{ij}, \lambda_{ij}): (\tU_j, \Gamma_j, \varphi_j)
\hookrightarrow (\tU_i, \Gamma_i, \varphi_i)$ we have 
${\psi_{ij}}^* \tf_i = \tf_j$.  Naturally each $\tf_i$ defines a
continuous map $f_i : \tU_i/\Gamma_i = U_i \to \R$.   Thanks to
the compatibility conditions above, these maps glue together to define a
continuous map from the topological space underlining $M$ to $\R$.  By
abuse of notation we may write $f:M \to \R$.  Similarly, a {\bf
differential $k$-form} $\sigma$ on the orbifold $M$ is a collection of
$\Gamma_i$-invariant $k$-forms $\tsigma_i$ on $\tU_i$ such that for
any injection $(\psi_{ij}, \lambda_{ij}): (\tU_j, \Gamma_j, \varphi_j)
\hookrightarrow (\tU_i, \Gamma_i, \varphi_i)$ we have 
${\psi_{ij}}^* \tsigma_i = \tsigma_j$.  We denote the collection of
all differential forms on $M$ by $\Omega^* (M)$.  Note that $\Omega^*
(M)$ has a well-defined exterior multiplication (since exterior
multiplication behaves well under pull-blacks) and exterior
differentiation $d : \Omega^* (M) \to \Omega ^{*+1} (M)$ (for the same
reason).
In the same manner one defines vector fields, Riemannian metrics,
quadratic forms and other differential geometric objects on orbifolds.
 
Finally we define maps of orbifolds.  The reader should be aware that
there are several notions of maps orbifolds.  For example Satake in
his two papers gave two inequivalent definitions.  The one we give
below is the simplest; it is not the best.  It will, however, suffice
for our purposes.  See Ruan's survey \cite{R} for the modern point of
view.  Let $(M, \{ \tU_i, \Gamma_i, \varphi_i\})$ and $(N, \{ \tV_j,
\Delta_j, \phi_j\})$ be two orbifolds and let $F: M \to N$ be a
continuous map of underlying topological spaces.  The map $F$ is a
(smooth) {\bf map of orbifolds} if for every point $x\in M$ there are charts
$(\tU_i, \Gamma_i, \varphi_i)$, $(\tV_j, \Delta_j, \phi_j)$ with $x
\in \varphi (\tU_i)$, $F (\varphi_i (\tU_i)) \subset \phi_j (\tV_j)$ 
and a $C^\infty $ map $\tilde{F}_{ji} : \tU_i \to \tV_j$ such that 
$$
\phi_j \circ \tilde{F}_{ji} = F\circ \varphi_i.
$$

The reason for the appearance of orbifolds in these notes is the
symplectic reduction theorem of Marsden, Weinstein and Meyer \cite{MW, Me}:
\begin{theorem}
Let $(M, \omega)$ be a symplectic manifold with a Hamiltonian action
of a Lie group $G$.  Let $\Phi : M\to \fg^*$ denote a corresponding
moment map.  Suppose $\eta \in \fg^*$ is a regular value of $\Phi$ and
suppose that the action of the isotropy group $G_\eta$ of $\eta$ on
$\Phi\inv (\eta)$ is proper.  Then the action of $G_\eta$ on $\Phi\inv
(\eta)$ is locally free and the quotient $M_\eta : = \Phi \inv
(\eta)/G_\eta$ is naturally a symplectic orbifold.
\end{theorem}
We will not discuss the proof of this well known theorem.  Instead let
me briefly explain why quotients by locally free actions are
orbifolds.  The reason is the slice theorem: if the action of a Lie
group $G$ is locally free and proper on a manifold $Z$, then for any
point $z\in Z$ there is a slice $S_z$ for the action of $G$ and the
isotropy group $G_z$ is finite.  The quotient $S_z/G_z$ is
homeomorphic to a neighborhood of the orbit $G\cdot z$ in the orbit
space $Z/G$.  The triple $(S_z, G_z, \varphi: S_z \to S_z/G_z
\hookrightarrow Z/G)$ is a uniformizing chart of the orbifold $Z/G$.

Let us now prove Theorem~\ref{BGthm}.
\begin{proof}[Proof of Theorem~\ref{BGthm}]
Since $B$ is compact, the image $\Psi_\alpha (B)$ is compact.
Therefore the set of vectors $X' \in \fg$, such that the function
$\langle \Psi_\alpha, X'\rangle$ is strictly positive on $B$, is open.
Hence we may assume that $X$ lies in the integral lattice $\Z_G :=
\ker (\exp: \fg \to G)$ of the torus $G$.  Let $H = \{\exp tX \mid
t\in \R\}$ be the corresponding circle subgroup of $G$.

Let $f(x) = 1/(\langle \Psi_\alpha (x) ,X\rangle)$ and let $\alpha' =
f \alpha$.  The form $\alpha '$ is another $G$-invariant contact form
with $\ker \alpha' = \xi$.  The moment map $\Psi_{\alpha'}$ defined by
$\alpha'$ satisfies $\Psi_{\alpha'} = f \Psi_\alpha$.  Therefore
$\langle \Psi_{\alpha'} (x), X\rangle = 1$ for all $x\in B$.
 
Since the function $\langle \Psi_\alpha , X \rangle $ is nowhere zero,
the action of $H$ on $B$ is locally free.  Consequently the induced
action of $H$ on the symplectization $(N, \omega) = (B\times \R,
d(e^t\alpha'))$ is locally free as well.  Hence any $a\in \R $ is a
regular value of the $X$-component $\langle \Phi, X\rangle$ of the
moment map $\Phi$ for the action of $G$ on the symplectization $(N,
\omega)$.  Note that $\Phi (x, t) = -e^t \Psi_{\alpha'}(x)$. The manifold 
$B\times \{0\}$ is the $-1$ level set of $\langle \Phi, X\rangle$.
Therefore by the reduction theorem $M:= (\langle \Phi, X\rangle )\inv
(-1)/H \simeq B/H$ is a (compact connected) symplectic orbifold.  The
action of $G$ on $(\langle \Phi, X\rangle )\inv$ descends to an
effective Hamiltonian action of $G/H$ on $M$.  A dimension count shows
that the effective action of $G/H$ on $M$ is completely integrable,
i.e., $M$ is a symplectic toric orbifold.
\end{proof}
 
The rest of the section is a proof that real odd dimensional cohomology
vanishes for symplectic toric orbifolds.   It uses Morse theory.

\section{Morse theory on orbifolds}

Let us start by briefly reviewing the fundamental results of Morse
theory on manifolds.  Let $f: M \to \R$ be a smooth function.  A {\bf
critical point} of $f$ is a point $p$ where the differential $df$ is
0. The image $f(p)$ of a critical point $p$ is a {\bf critical value}
of $f$.  A critical point $p$ of $f$ is {\bf nondegenerate} if the
Hessian $d^2 f _p$ is a nondegenerate quadratic form (in local
coordinates $d^2f_p$ is the matrix of second order partials
$\left(\frac{\partial^2 f}{\partial x_i \partial x_j}\right)$; thus
$d^2 f_p (x) = \sum _{ij}\frac{\partial^2 f}{\partial x_i \partial
x_j} x_i x_j$).  The {\bf index} of a nondegenerate critical point $p$
is the number of negative eigenvalues of the Hessian $d^2 f_p$
(counted with multiplicities.)  The two and a half fundamental results
of Morse theory are the two theorems and the lemma below.
\begin{theorem}\labell{MthmA}
Let $f$ be a smooth function on a manifold $M$, and $M_a$ the set
$f\inv(\,(-\infty, a]\,)$. If $f\inv ([a,b])$ is compact and contains no
critical points then $M_a$ is homotopy equivalent to $M_b$.
\end{theorem}

\begin{proof}[Sketch of proof] Fix a Riemannian metric on $M$.  Since $f$ has 
no critical points in $f\inv ([a, b])$, the unit vector field $X = -
\nabla f /||\nabla f||$ is well-defined.  Extend $X$ to all of $M$.
The flow of $X$ gives a retraction of $M_b$ onto $M_a$.
\end{proof}

\begin{theorem}\labell{MthmB}
Let $f$ be a smooth function on a manifold $M$.  Suppose $f\inv
([a,b])$ is compact and contains exactly one nondegenerate critical
point $p$ in its interior.  Then $M_b$ has the homotopy type of $M_a$
with a $\lambda$-dimensional disk $D^\lambda$ attached along the
boundary $\partial D^\lambda = S^{\lambda -1}$ where $\lambda$ is the
index of $p$.
\end{theorem}
We omit the proof of the theorem which is well known noting only that
the key ingredient of the proof is
\begin{lemma}[Morse Lemma] \labell{morse lemma}
Let $p$ be a nondegenerate critical point of index $\lambda$ of a
function $f: M \to \R$.  There is a neighborhood $U$ of $p$ in $M$ and
an open embedding $\varphi: U \to M$ such that $f\circ \varphi (x) =
f(p) + d^2f_p (x)$ for all $x\in U$.  There is, there is a change of
coordinates near $p$ so that in new coordinates $f$ is a quadratic
form (up to a constant).

Moreover, if a compact Lie group $G$ acts on $M$ fixing $p$ and
preserving $f$, we may arrange for $U$ to be $G$-invariant and for
$\varphi$ to be $G$-equivariant.  Note that in this case the Hessian
$d^2f_p (x)$ is a $G$-invariant quadratic form.
\end{lemma}

Suppose now that $f$ is a smooth function on an {\em orbifold} $M$.
Then the 1-form $df$ still makes sense.  We define $p$ to be a {\bf
critical point} of $f$ if $df_p =0$.  A critical point $p$ is {\bf
nondegenerate} if for any uniformizing chart $(\tU, \Gamma, \varphi)$
with $p \in U = \varphi (\tU) $, a point $\tilde{p} \in \varphi \inv
(p)$ is a nondegenerate critical point of $\tilde{f} = f\circ
\varphi$.  The {\bf index} of $p$ is the index of the
$\Gamma$-invariant quadratic form $d^2 \tilde{f}_{\tilde {p}}$.

It is not hard to believe that Theorem~\ref{MthmA} holds for orbifolds
with no changes and that essentially the same proof still works.
Theorem~\ref{MthmB} requires a small modification, see \cite{LT}.

\begin{theorem}\labell{MthmB'}
 Let $f$ be a smooth function on an orbifold $M$.  Suppose $f\inv
 ([a,b])$ is compact and contains exactly one nondegenerate critical
 point $p$ in its interior.  Then $M_b$ has the homotopy type of $M_a$
 with the quotient $D^\lambda/\Gamma$ attached along the boundary
 $(\partial D^\lambda)/\Gamma = S^{\lambda -1}/\Gamma$ where $\lambda$
 is the index of $p$ and $\Gamma$ is structure group of $p$. Here
 again $D^\lambda$ and $S^{\lambda -1}$ denote the disk of dimension
 $\lambda$ and the sphere of dimension $\lambda -1$ respectively.
\end{theorem}
\begin{corollary}\label{CorA}
Let $f: M \to \R$, $p$, $\lambda$ and $\Gamma$ be as above.  Then $H^*
(M_b, M_a; \R) = H^* (D^\lambda/\Gamma , S^{\lambda -1}/\Gamma ; \R)$.
Hence, if $\Gamma$ is orientation preserving, $H^q (M_b, M_a; \R) =
\R$ for $q = \lambda$ and 0 otherwise.
\end{corollary}
\begin{proof}
By excision $H^* (M_b, M_a; \R) = H^* (D^\lambda/\Gamma , S^{\lambda
-1}/\Gamma ; \R)$.  If $\Gamma$ is orientation preserving, then $H^q
(D^\lambda/\Gamma , S^{\lambda -1}/\Gamma ; \R) = \tilde{H}^q
(S^\lambda/\Gamma ; \R) = \R$ for $q = \lambda$ and 0 otherwise. Here
$\tilde{H}^q$ denotes the reduced cohomology.
\end{proof}
\noindent
As a consequence of the corollary above we get 
\begin{corollary}\label{CorB}
 Let $f$ be a smooth function on a compact orbifold $M$.  Suppose all
 indices of $f$ are even.  Then $H^q (M, \R) = 0$ for all odd indices
 $q$.
\end{corollary}
\begin{proof}
The proof is inductive.  Suppose that $c$ is a critical value of $f$
and suppose we know that for all $a< c$ we have $H^q (M_a, \R) = 0$
for $q$ odd.  Assume for simplicity that $f\inv (c)$ contains only one
critical point (this keeps the notation more manageable) and that its
index is $2k$ for some integer $k$. Then using the long exact sequence
for the pair $(M_{c + \epsilon}, M_{c-\epsilon})$ (for some
sufficiently small $\epsilon >0$) and Corollary~\ref{CorA} we see that
$H^q (M_{c+\epsilon}, \R) = H^q (M_{c-\epsilon}, \R)$ for $q\not =
2k-1, 2k$, that $H^{2k-1} (M_{c+\epsilon}, \R)$ embeds in $H^{2k-1}
(M_{c-\epsilon}, \R) = 0$ and that $H^{2k} (M_{c+\epsilon}, \R) =
H^{2k} (M_{c-\epsilon}, \R)\oplus \R$.  The result follows.
\end{proof}

Now suppose $(M, \omega, \Phi:M \to \fg^*)$ is a compact symplectic
toric orbifold.  Just as for symplectic toric manifolds the image
$\Phi(M)$ is a simple polytope and the moment map sends the fixed
points $M^G$ in one-to-one fashion to vertices of $\Phi (M)$ (see \cite{LT}).
Therefore, for a generic vector $X\in \fg$ the function $f = \langle
\Phi, X \rangle$ takes distinct values at fixed points.  We will now argue 
that the critical points of $f$ are exactly the fixed points.  We will
then argue that $f$ is Morse and that all indices of $f$ are even.

For an action of a torus $G$ on a compact orbifold $M$ only finitely
many isotropy groups can occur.  This is a consequence of compactness
and existence of slices.  Therefore the set of subalgebras $\fg_x$,
which are Lie algebras of isotropy groups for the action of $G$ on
$M$, is finite.  Hence the set
$$
\cU = \fg \smallsetminus \bigcup_{x\in M} \{\fg_x \mid \fg_x \not = \fg\}
$$
is open and dense.  Now take  $X$ to be a vector in $\cU$.  Then 
\begin{equation*}
\begin{split}
0= df_x = d \langle \Phi, X\rangle _x = (\iota (X_M) \omega)_x 
& \Leftrightarrow X_M (x) = 0\\
	& \Leftrightarrow \exp tX \cdot x = x \text{ for all } t \\
	& \Leftrightarrow \exp tX \in G_x  \text{ for all } t\\
	& \Leftrightarrow X \in \fg_x .\\
\end{split}
\end{equation*}
Since $X\in \cU$ we see that $df_x = 0$ $\Leftrightarrow$ $\fg_x =
\fg$ $\Leftrightarrow$ $x$ is fixed by $G$.
It remains to check that $f$ is Morse and that all the indices of $f$
are even.

Let $x\in M^G$ be a fixed point and let $(\tU , \Gamma, \varphi)$ be a
uniformizing chart centered at $x$ (cf.\ Remark~\ref{centered}).  We
may assume that $U = \varphi (\tU)$ is $G$-invariant.  Denote the
symplectic form on $\tU$ by $\tilde{\omega}$.  The map $\tPhi = \Phi
\circ \varphi : \tU \to
\fg^*$ is a moment map for an action of a torus $\tilde{G}$ on $\tU$.
One can show arguing as in the proof of Lemma~\ref{Del-lemma} that
$\Gamma \subset \tilde{G}$ and that $G = \tilde{G}/\Gamma$.  In
particular $\tilde{G}$ has $\fg$ as its Lie algebra.  We now apply the
equivariant Darboux theorem to $\tU$, $\tilde{\omega}$ and
$\tilde{G}$.  We get a $\tilde{G}$ invariant neighborhood $\tV$ of
$\tilde{x} = \varphi\inv (x)$, an open neighborhood $\tV_0$ of 0 in
$T_{\tilde{x}} \tU$ and a $\tilde{G}$-equivariant diffeomorphism $\tau:
\tV_0 \to \tV$ such that $\tau^* \tilde{\omega} = \tilde{\omega}_0$
where $\tilde{\omega}_0$ is the constant coefficient form
$\tilde{\omega}_{\tilde{x}}$ on the vector space $T_{\tilde{x}} \tU$.

The action of $\tilde{G}$ on $T_{\tilde{x}} \tU$ is linear.  Hence the
corresponding moment map $\tPhi_0 :T_{\tilde{x}} \tU \to \fg^*$ is
quadratic. In fact by choosing a $\tilde{G}$-invariant complex
structure on $T_{\tilde{x}} \tU$ compatible with the symplectic form we
can identify $(T_{\tilde{x}} \tU, \tilde{\omega}_0)$ with $(\C^n,
\sqrt{-1} \sum dz_j \wedge d \bar{z}_j)$ so that the action of
$\tilde{G}$ is given by
$$
a \cdot (z_1, \ldots, z_n)=  (\chi_1 (a) z_1, \ldots, \chi_n (a)z_n) 
$$
for some characters $\chi_j : \tilde{G} \to S^1$.  Then 
$$
\tPhi_0 (z_1, \ldots, z_n) = \sum |z_j|^2 \nu_j
$$ where $\nu_j = d\chi_j$ are the corresponding weights.  Since
$\tau$ is a $\tilde{G}$-equivariant symplectomorphism, $\tPhi \circ
\tau = \tPhi_0 + c$ for some constant $c\in \fg^*$.  Hence
$$
f \circ \varphi \circ \tau = \langle \Phi, X\rangle  \circ \varphi \circ \tau 
= \langle \tPhi, X\rangle  \circ \tau = \langle \tPhi_0, X\rangle + c
\text{, i.e.,}
$$
$$
f \circ \varphi \circ \tau (z_1, \ldots, z_n) = \sum \nu_j (X) |z_j|^2 + c .
$$
By the choice of $X$, $\nu_j (X) \not = 0$ for any $j$.  Thus $f\circ
\varphi$ is Morse and its index at $\tilde{x}$ is twice the number of
weights $\nu_j$ with $\nu_j (X) < 0$.  Therefore, by definition, $f$
is a Morse function on $M$ and all indices of $f$ are even.

\appendix
\chapter{Hypersurfaces of contact type}\label{appA}
Contact manifolds often arise as codimension 1 submanifolds of
symplectic manifolds, i.e., as hypersurfaces.

\begin{definition}
Let $(M, \omega)$ be a symplectic manifold.  A hypersurface $\Sigma $
of $M$ is of {\bf contact type} if there is a neighborhood $U$ of
$\Sigma$ in $M$ and a vector field $X$ on $U$ such that
\begin{enumerate}
\item $T_m M = T_m \Sigma \oplus \R X(m)$ for any point $m\in \Sigma$, 
i.e., the vector field $X$ is nowhere tangent to $\Sigma$;
\item  the flow of $X$ expands the symplectic form $\omega$ exponentially, 
i.e., $L_X \omega = \omega$, where as usual $L_X$ denotes the Lie
derivative with respect to $X$.
\end{enumerate}
The vector field $X$ with above properties is often called a {\bf
Liouville vector field}.
\end{definition}

We now prove that hypersurfaces of contact type are indeed contact manifolds.
\begin{proposition}
Let $\Sigma$ be a hypersurface of contact type in a symplectic
manifold $(M,\omega)$ and let $X$ be a Liouville vector field defined
on a neighborhood $U$ of $\Sigma$.  The 1-form $\alpha =
(\iota(X)\omega) |_\Sigma$ is contact.
\end{proposition}
\begin{proof}
Note first that $d (\iota (X) \omega) = d (\iota (X) \omega) + \iota
(X) d \omega = L_X \omega = \omega$.  Hence $d\alpha = d
(\iota(X)\omega) |_\Sigma = \omega |_\Sigma$. Since $\omega$ is
symplectic and $\Sigma$ is of codimension 1 in $M$, the form $\omega_m
|_{T_m \Sigma}$ has a 1-dimensional kernel (for any point $m\in
\Sigma$). Thus there is a vector $Y_m \in T_m \Sigma$ such that
$\omega _m (Y_m, v) = 0$ for any $v\in T_m \Sigma$. Since $\omega$ is
symplectic and since by assumption $T_m M = T_m \Sigma \oplus \R X(m)$
for all $m\in \Sigma$ we have
$$
0 \not = \omega_m (Y_m, X(m)) = -(\iota (X) \omega)_m (Y_m) 
= - \alpha _m (Y_m). 
$$
Hence $\alpha _m \not = 0$ for all $m\in \Sigma$ and consequently $\xi
= \ker \alpha$ is a codimension 1 distribution on $\Sigma$. It remains to 
show that for any $m\in \Sigma$ 
$$
 d\alpha _m |_{\xi _m} = \omega_m |_{\xi _m}
$$
is non-degenerate.  On the other hand $d \alpha_m = \omega _m |_{T_m
\Sigma}$ and $\xi _m = \{ v\in T_m \Sigma \mid \omega_m (v, X(m)) =
0\}$.  Since $\omega _m (v, Y_m) = 0$ for any $v\in T_m \Sigma$, the subspace 
$\xi_m $ of $T_m M$ lies in the symplectic perpendicular to the 2-plane 
$\text{span}_\R \{Y_m, X(m)\}$ in $(T_m M, \omega_m)$:
$$
\xi _m \subset (\text{span}_\R \{Y_m, X(m)\})^\omega .
$$ 

Since $\dim \xi_m = \dim M - 2 = \dim (\text{span}_\R \{Y_m,
X(m)\})^\omega$, we have equality: $\xi _m = (\text{span}_\R \{Y_m,
X(m)\})^\omega$.  Since $(\text{span}_\R \{Y_m, X(m)\})$ is a
symplectic subspace, its symplectic perpendicular $\xi_m$ is
symplectic as well.
\end{proof}
\begin{example}
Let $(Q, g)$ be a Riemannian manifold.  Let $g^*$ denote the dual
metric on $T^*Q$.  The co-sphere bundle $S^*Q$ is the set of covectors
in $T^*Q$ of length 1: $S^*Q = \{(q,p) \in T^*Q \mid q\in Q, \, p\in
T^*_q Q, \, g^*_q (p,p) = 1\} $.  It is a hypersurface of contact type
in $T^*Q$ relative to the standard symplectic structure.  The
Liouville vector field is the generator of dilations.  In local
coordinates $X(q, p) = \sum p_j \frac{\partial}{\partial q_j}$.
\end{example}
\begin{exercise}
A codimension 1 hypersurface $\Sigma \subset \C^n$ is {\bf star-shaped
about the origin} if for any nonzero vector $v\in \C^n$ the ray $\{ tv
\mid t \in (0, \infty)\}$ intersects $\Sigma$ transversely in exactly one 
point.  In particular $\Sigma$ is the image of an embedding $\iota :
S^{2n -1} \hookrightarrow \C^n$ of the $(2n-1)$-dimensional sphere.

Show that any star-shaped hypersurface is a hypersurface of contact
type in $(\C^n, \omega = \sqrt{-1} \sum dz_j \wedge d \bar{z}_j)$.
Show that any two star-shaped hypersurfaces are isomorphic as contact
manifolds.  The contact structure in question is called the {\bf
standard contact structure on $S^{2n-1}$}.

Prove a converse: given a contact form $\alpha $ on $S^{2n-1}$
defining the standard contact structure there is an embedding $\imath :
S^{2n-1} \to \C^n$ such that $\imath^* (\iota(X) \omega) = \alpha$
where $X$ is the radial vector field on $\C^n$: $X (z) = \frac{1}{2}
\sum \left( z_j \frac{\partial }{\partial z_j} + \bar{z}_j
\frac{\partial }{\partial {z}_j}\right)$.

\end{exercise}

\end{document}